\documentclass[psamsfonts]{amsart}
\usepackage{amssymb}
\usepackage{graphicx,color}
\usepackage{amssymb,amsfonts,amsthm}
\usepackage[all,arc]{xy}
\usepackage{enumerate}
\usepackage{mathrsfs}

\newtheorem{thm}{Theorem}[section]

\newtheorem{lemma}[thm]{Lemma}
\newtheorem{prop}[thm]{Proposition}

\theoremstyle{definition}

\theoremstyle{remark}
\newtheorem{rem}[thm]{Remark}
\makeatletter
\let\c@equation\c@thm
\makeatother
\numberwithin{equation}{section}

\title[Asymptotic Plateau problem]{Asymptotic Plateau problem for $3$-convex hypersurface in $\mathbb{H}^5$}

\author{Zhenan Sui}

\address{Institute for Advanced Study in Mathematics of HIT, Harbin Institute of Technology, Harbin, China}
\email{20170045@hit.edu.cn}

\begin{document}

\begin{abstract}
We prove the existence of a smooth complete $3$-convex hypersurface which satisfies prescribed curvature equation $\prod\limits_{i = 1}^n (H - \kappa_i) = \big( (n - 1) \sigma \big)^n$ for $n = 4$ and has prescribed asymptotic boundary $\Gamma$ at the infinity of hyperbolic space of dimension 5, where $\sigma \in (0, 1)$ is a constant and $\Gamma$ is assumed to have nonnegative mean curvature. We introduce Lagrange multiplier method to compute the extreme value of the concavity of $f (\kappa) = \frac{1}{n - 1} \Big( \prod\limits_{i = 1}^n (H - \kappa_i) \Big)^{\frac{1}{n}}$ during uniform global curvature estimate.
\end{abstract}

\subjclass[2010]{Primary 53C21; Secondary 35J65, 58J32}

\maketitle


\section {\large Introduction}

\vspace{4mm}

This paper is devoted to uniform global curvature estimate for asymptotic Plateau problem in hyperbolic space, which is an open problem in Guan and Spruck \cite{GS10}. Combined with the other a priori estimates established in \cite{GS10}, one can prove the existence of classical solutions to asymptotic Plateau problem.

We shall first describe asymptotic Plateau problem in hyperbolic space: in the ambient space
\[\mathbb{H}^{n+1} = \big\{ (x, x_{n+1}) \in \mathbb{R}^{n+1} \big\vert x_{n + 1} > 0 \big\},
\quad d s^2 = \frac{1}{x_{n+1}^2} \sum_{i = 1}^{n + 1} d x_i^2, \]
we want to find a smooth complete hypersurface $\Sigma$ which satisfies the equation
\begin{equation} \label{eq1-1}
f \big( \kappa[ \Sigma ] \big) = \sigma
\end{equation}
and is subject to the asymptotic boundary condition
\begin{equation} \label{eq1-2}
\partial \Sigma = \Gamma,
\end{equation}
where $\kappa[ \Sigma ] = (\kappa_1, \cdots, \kappa_n)$ are the principal curvatures of $\Sigma$, $f$ is a smooth symmetric function of $n$ variables, $\sigma \in (0, 1)$ is a prescribed constant and $\Gamma$ is a disjoint collection of smooth closed embedded $(n - 1)$ dimensional submanifolds in
\[ \partial_{\infty} \mathbb{H}^{n + 1} = \mathbb{R}^n \times \{ 0 \} \equiv \mathbb{R}^n.\]

This problem was initially studied by Anderson \cite{Anderson1982, Anderson1983}, and Hardt and Lin \cite{HL1987} for area-minimizing varieties by geometric measure theory, which were generalized by Tonegawa \cite{Tonegawa1996} to hypersurfaces of constant mean curvature. Later, PDE methods were introduced by Lin \cite{Lin1989} for minimal hypersurfaces, followed by Nelli and Spruck \cite{NS1996} and Guan and Spruck \cite{GS00} for hypersurfaces of constant mean curvature. The Gauss curvature case was first studied by Labourie \cite{Labourie1991} in $\mathbb{H}^3$, and then in $\mathbb{H}^{n + 1}$ by Rosenberg and Spruck \cite{RS1994}. For general curvature function $f$ defined on
the positive cone $\Gamma_n$,
asymptotic Plateau problem was completely solved by Guan, Spruck, Szapiel and Xiao \cite{GSS09, GS11, GSX14} under several fundamental structure conditions. Recently, Hong, Li and Zhang \cite{HLZ2025} introduced an alternative method for constructing strictly locally convex solutions to \eqref{eq1-1}-\eqref{eq1-2}.

For general $f$ defined on general cone $K$, in \cite{GS10}, Guan and Spruck proved:
\begin{thm}\textsuperscript{\cite{GS10}} \label{GuanSpruckTheorem}
Let $\Gamma = \partial \Omega \times \{ 0 \} \subset \mathbb{R}^{n + 1}$ where $\Omega$ is a bounded smooth domain in $\mathbb{R}^n$. Suppose that the Euclidean mean curvature $\mathcal{H}_{\partial \Omega}$ is nonnegative and $\sigma \in (0, 1)$ satisfies $\sigma > \sigma_0$, where $\sigma_0$ is the unique zero in $(0, 1)$ of
\[ \phi(a) := \frac{8}{3} a + \frac{22}{27} a^3 - \frac{5}{27} (a^2 + 3)^{3/2}. \]
(Numerical calculations show $0.3703 < \sigma_0 < 0.3704$.)  Under conditions \eqref{eq2-16}--\eqref{eq2-21},
there exists a complete hypersurface $\Sigma$ in $\mathbb{H}^{n + 1}$ satisfying \eqref{eq1-1}--\eqref{eq1-2} with uniformly bounded principal curvatures
\begin{equation} \label{eq1-7}
\big| \kappa [ \Sigma ] \big| \leq C \quad \text{ on } \Sigma.
\end{equation}
Moreover, $\Sigma$ is the graph of a unique admissible solution $u \in C^{\infty} ( \Omega ) \cap C^1( \overline{\Omega} )$ of the Dirichlet problem associated to \eqref{eq2-14}
\begin{equation*}
\left\{ \begin{aligned}
G (D^2 u, D u, u) = & f \big( \kappa [ u ] \big) = \sigma, \quad u > 0 \quad \text{ in }  \Omega \subset \mathbb{R}^n \\
u = & 0 \quad \text{ on } \partial\Omega.
\end{aligned} \right.
\end{equation*}
Furthermore, $u^2 \in C^{\infty} (\Omega) \cap C^{1, 1} (\overline{\Omega})$ and
\[ \sqrt{1 + |D u|^2} \leq \frac{1}{\sigma}, \quad u |D^2 u| \leq C \quad \text{ in } \Omega, \]
\[ \sqrt{1 + |D u|^2} = \frac{1}{\sigma} \quad \text{ on } \partial \Omega. \]
\end{thm}

Our aim is to derive uniform global curvature estimate without assumption $\sigma > \sigma_0$. All the other a priori estimates up to second order have been derived in \cite{GS10} without assumption $\sigma > \sigma_0$. Throughout this paper, we always assume that $\Sigma$ can be represented as a vertical graph $u$ over the smooth bounded domain $\Omega \subset \mathbb{R}^n$ surrounded by $\partial \Omega = \Gamma$.  We also assume $\Gamma$ to be mean convex as \cite{GS10}.

The essential difficulty for solving asymptotic Plateau problem \eqref{eq1-1}--\eqref{eq1-2} in hyperbolic space is the singularity of the equation at $\Gamma$. Because the prescribed curvature is a constant, and the asymptotic boundary $\Gamma$ is assumed to be mean convex when $K \supsetneq \Gamma_n$, the key ingredients for the proof are:  (1) prove the existence of solution $u^{\epsilon}$ to the approximate Dirichlet problem
\begin{equation} \label{eq1-4}
\left\{ \begin{aligned}
f( \kappa [ u ] ) = & \sigma \quad \text{ in } & \Omega, \\
u = & \epsilon \quad \text{ on } & \Gamma,
\end{aligned} \right.
\end{equation}
(2) give uniform a priori estimates up to second order for solution sequence $u^{\epsilon}$ to guarantee a convergent subsequence whose limit is a solution to \eqref{eq1-1}--\eqref{eq1-2}. Here ``uniform'' means the estimates should be independent of $\epsilon$; certainly, the estimates may depend on $\sigma$.

Progress towards this direction, that is, to derive uniform global curvature estimate for any prescribed $\sigma \in (0, 1)$ is first made by Lu \cite{Lu2023} for $f = \sigma_{n - 1}$,
where
\[\sigma_k (\kappa) = \sum_{1 \leq i_1 < \ldots < i_k \leq n} \kappa_{i_1} \cdots \kappa_{i_k} \]
is the $k$th elementary symmetric function which is defined on $k$-th G\r arding's cone
\[\Gamma_k = \big\{ \kappa \in \mathbb{R}^n \big\vert \sigma_j (\kappa) > 0, \, j = 1, \ldots, k \big\}. \]
A key observation is: one can use Guan, Ren and Wang's inequalities \cite{GRW2015, RW2019, RW2023} to tackle third order terms delicately; in other words, concavity can be well exploited to balance the difficult terms arising from the hyperbolic space. In this spirit, Wang solved the case when $f = \frac{\sigma_k}{\sigma_{k - 1}}$ in \cite{WangB-MRL}, $f = \sigma_2$ in \cite{WangB-Adv} and $f = \sigma_{n - 2}$ in \cite{WangB-arXiv}. In particular, Shankar and Yuan's almost Jacobi inequality \cite{SY2025} were fully utilized in \cite{WangB-Adv}. These successes reflect that concavity plays a crucial role in uniform global curvature estimate.

In this paper, we continue our research on the following prescribed curvature equation:
\begin{equation} \label{eq1-8}
\prod\limits_{i = 1}^n (H - \kappa_i) = \big( (n - 1) \sigma \big)^n.
\end{equation}
Here
\[ H = \sum\limits_{i = 1}^n \kappa_i \]
is the mean curvature of $\Sigma$. Admissible functions subject to equation \eqref{eq1-8} are called $(n - 1)$-convex hypersurfaces, which were studied intensively by Sha \cite{Sha1986, Sha1987}, Wu \cite{Wu1987} and Harvey and Lawson \cite{HarveyLawson2013}. In the complex setting, equation \eqref{eq1-8} is related to the Gauduchon conjecture \cite{Gauduchon1984}, which was solved by Sz\'ekelyhidi, Tosatti and Weinkove \cite{STW2017}.

Equation \eqref{eq1-8} happens to be the intersection of the following two classes of curvature functions.
\begin{description}
  \item[Class 1\textsuperscript{\cite{CJ2021}}]
\begin{equation} \label{eq1-5}
f ( \kappa ) = \frac{1}{(C_n^k)^{\frac{1}{k}} (n - 1)} \sigma_k^{\frac{1}{k}}(\lambda), \qquad \lambda \in \Gamma_k,  \qquad 1 \leq k \leq n,
\end{equation}
where
\[ \lambda_i = \sum\limits_{j \neq i} \kappa_j, \quad i = 1, \ldots, n, \]
  \item[Class 2\textsuperscript{\cite{Dinew2023, Dong2023}}]
\begin{equation} \label{eq1-6}
f ( \kappa ) = \frac{1}{k} P_k( \kappa )^{\frac{1}{C_n^k}},
\end{equation}
where
\[ P_k (\kappa) = \prod\limits_{i_1 < \cdots < i_k} (\kappa_{i_1} + \cdots + \kappa_{i_k}), \qquad 1 \leq k \leq n \]
is defined on the cone
\[ \mathcal{P}_k = \big\{ \kappa \in \mathbb{R}^n \big| \kappa_{i_1} + \cdots + \kappa_{i_k} > 0, \ \ \forall \, 1 \leq i_1 < \cdots < i_k \leq n  \big\}. \]
\end{description}
Both classes satisfy structure condition \eqref{eq2-16}--\eqref{eq2-21} as in Guan and Spruck \cite{GS10} (one may see \cite{CJ2021, Dinew2023, Dong2023} for the verifications). Admissible hypersurface $\Sigma$ associated to \eqref{eq1-6} (that is, $\kappa[\Sigma] \in \mathcal{P}_k$) is called $k$-convex hypersurface.
Recent work on equation \eqref{eq1-5} and \eqref{eq1-6} can be found in literature
\cite{CDH2023, Chen-Tu-Xiang, Chen-Tu-Xiang2021, Chen-Tu-Xiang2023, Dong2024, Jiao-Sun, Lu-Zhong, Mei2023, Mei-Zhu, Sheng-Xue, Wang2025, Zhou2022, Zhou2024}.

A significant difference on uniform global curvature estimate for asymptotic Plateau problem in hyperbolic space is: we can no longer use any zero order term in test function, for otherwise the estimate would depend on $\epsilon$, which is not uniform. This restriction may bring in huge difficulties in the estimation. In Chen, Sui and Sun \cite{CSS2025}, we found that
for curvature function \eqref{eq1-5} with $k < n$, uniform global curvature estimate can be easily obtained by applying the special property $f_1 \geq C \sum f_i$ (assuming $\kappa_1 \geq \cdots \geq \kappa_n$) due to Chen, Dong and Han \cite{CDH2023};
for $k = n$, \eqref{eq1-5} no longer possesses this special property, and uniform global curvature estimate becomes extremely difficult. Because the experiences in all previous literatures can not be applied, in \cite{CSS2025}, we computed the concavity accurately by solving six equations to tackle the special case $k = n = 3$, which is the same with $k = 2$ and $n = 3$ in \eqref{eq1-6}.

Following the route of \cite{CSS2025}, in this paper, we investigate asymptotic Plateau problem associated to equation \eqref{eq1-8} with $n = 4$. There are two huge challenges:
\begin{description}
  \item[Question 1] how to compute concavity accurately for fully nonlinear operators?
  \item[Question 2] how to analyze a complicated quantity arising from concavity and how to give a suitable classification to realize the estimation?
\end{description}

In this paper, Question 1 is completely solved for equation \eqref{eq1-8} by Lagrange multiplier method. To be more concrete, we find the extreme value of the quadratic form involving third order terms under the constraint \eqref{eq4-23} and \eqref{eq4-6} for each fixed $i$. This result is summarized in Theorem \ref{OptimalConcavity}. Besides the operator associated to \eqref{eq1-8}, Lagrange multiplier method also works for $\sigma_2$ and $\frac{\sigma_2}{\sigma_1}$ operators, which can lead to interior second order estimates for the associated equations (see Jiao and Sui \cite{JS2026}). In \cite{JS2026}, we also found that Lagrange multiplier method can recover Shankar and Yuan's almost Jacobi inequality in a different form. One may find it interesting to explore more application of Lagrange multiplier method.

For Question 2, as the quantity appeared in Theorem \ref{OptimalConcavity} is complicated, yet we believe that this extreme value of concavity should be accurately calculated so as to play its effect, therefore we confine ourselves to $n = 4$ and use Mathematica to obtain its exact value. Indeed, we have used Mathematica to conduct most of the calculation since Remark \ref{n=3}. For $n = 4$, we do the estimation for each $i$ layer with $i = 1, \ldots, 4$. We observe that each layer is different from the others, and we need to divide the estimation into more cases for discussion compared to dimension $n = 3$, which convinces us that the $n = 4$ case is richer than $n = 3$. This method can be generalized to general $n$ for equation \eqref{eq1-8}. Despite the vast computation involved, one can see the beautiful bridge between the geometric equation and the hidden algebraic structure.

Now we state the main result of this paper.
\begin{thm} \label{Theorem2}
Let $\Gamma = \partial \Omega \times \{ 0 \} \subset \mathbb{R}^{5}$ where $\Omega$ is a smooth bounded domain in $\mathbb{R}^4$. Suppose that $\sigma \in (0, 1)$ is a constant and the Euclidean mean curvature of $\partial \Omega$ is nonnegative. Then there exists a complete hypersurface $\Sigma$ in $\mathbb{H}^{5}$ satisfying \eqref{eq1-8}--\eqref{eq1-2} with uniformly bounded principal curvatures
\begin{equation} \label{eq1-9}
\big| \kappa [ \Sigma ] \big| \leq C \quad \text{ on } \Sigma.
\end{equation}
The rest of the conclusions are the same with Theorem \ref{GuanSpruckTheorem}.
\end{thm}
Consistent with our observation in the $n = 3$ case, for $n = 4$, we have to make each step optimal in order to make the estimation possible. We also observe that more parameters are involved for $n = 4$ and the choice of the range of parameters are different from $n = 3$ but is consistent with $n = 3$. It is an interesting question to investigate uniform global curvature estimate for equation \eqref{eq1-6}.

The rest of this paper is organized as follows: in section 2, we collect some calculations and facts for asymptotic Plateau problem in hyperbolic space. In section 3, we begin uniform global curvature estimate for \eqref{eq1-8}. In Section 4, we use Lagrange multiplier method to compute the extreme value of concavity accurately. From section 5 to section 8, we do the estimation for each $i$ layer. In section 9, we make the choice of the parameters and finish the proof of Theorem \ref{Theorem2}.

\vspace{4mm}

\section{Preliminaries}

\vspace{4mm}
In this section, we present some known calculations and facts for asymptotic Plateau problem in hyperbolic space. The proof can be found in \cite{GSS09, GS10, GS11, GSX14, CSS2025}.

\vspace{2mm}

\subsection{Geometric quantities on vertical graphs}~

\vspace{2mm}

Throughout this paper, we assume that $\Sigma$ can be represented as a vertical graph
\[ \Sigma = \Big\{ \big( x, u(x) \big) \big| x \in \Omega \Big\},  \]
where $\Omega$ is a smooth bounded domain in $\mathbb{R}^n \equiv \partial_{\infty} \mathbb{H}^{n + 1}$.

Denoting
\[ \partial_{I} = \frac{\partial}{\partial x_{I}}, \quad I = 1, \ldots, n + 1 , \]
then the coordinate tangent vector fields on $\Sigma$ can be written as
\[  \partial_i + u_i \partial_{n + 1}, \quad i = 1, \ldots, n. \]

If $\Sigma$ is viewed as a graph in Euclidean space $\mathbb{R}^{n + 1}$, then its upward unit normal vector field,  metric, inverse of the metric and second fundamental form are given by
\[\nu = \frac{( - D u, 1 )}{w}, \quad w = \sqrt{ 1 + |D u |^2},  \]

\[ \tilde{g}_{ij} = \delta_{ij} + u_i u_j, \]

\[ \tilde{g}^{ij} = \delta_{ij} - \frac{u_i u_j}{w^2},  \]

\[ \tilde{h}_{ij} = \frac{u_{ij}}{w}. \]
Therefore, by definition, the Euclidean principal curvatures $\tilde{\kappa} = (\tilde{\kappa}_1, \cdots, \tilde{\kappa}_n)$ of $\Sigma$ are the eigenvalues of
\[ \tilde{A} = \{ \tilde{a}_{ij} \}, \]
where
\[ \tilde{a}_{ij} = \frac{1}{w} \gamma^{ik} u_{kl} \gamma^{lj}, \]

\[ \gamma^{ik} = \delta_{ik} - \frac{u_i u_k}{w ( 1 + w )}, \]

\[ \gamma_{ik} = \delta_{ik} + \frac{u_i u_k}{1 + w}. \]
One can see that these quantities satisfy
\[ \gamma^{ik} \gamma_{kj} = \delta_{ij} \]
and
\[ \gamma_{ik} \gamma_{kj} = \tilde{g}_{ij}.  \]

When the ambient space is $\mathbb{H}^{n + 1}$, the unit upward normal vector field, metric, second fundamental form of $\Sigma$ become respectively
\[ {\bf n} = u \nu, \]

\[ g_{ij} = \frac{ \delta_{ij} + u_i u_j }{u^2}, \]

\[ h_{ij} = \frac{ \delta_{ij} + u_i u_j + u u_{ij} }{u^2 w} . \]
By definition, the hyperbolic principal curvatures $\kappa = (\kappa_1, \cdots, \kappa_n)$ are the eigenvalues of
\[ A = \{ a_{ij} \}, \]
where
\begin{equation} \label{eq2-13}
\begin{aligned}
a_{ij} = & u^2 \gamma^{ik} h_{kl} \gamma^{lj}
       =  \frac{\gamma^{ik} \big( \delta_{kl} + u_k u_l + u u_{kl} \big) \gamma^{lj}}{w}
       =  \frac{\delta_{ij} + u \gamma^{ik} u_{kl} \gamma^{lj}}{w} .
\end{aligned}
\end{equation}

Denoting $\nu^{n + 1} = \nu \cdot \partial_{n + 1}$, where $\cdot$ is the inner product in $\mathbb{R}^{n + 1}$. By the above discussion, we obtain the relations
\begin{equation} \label{eq2-11}
h_{ij} = \frac{1}{u} \tilde{h}_{ij} + \frac{\nu^{n+1}}{u^2} \tilde{g}_{ij}
\end{equation}
and
\begin{equation} \label{eq2-12}
\kappa_i =  u \tilde{\kappa_i} + \nu^{n+1},  \quad i = 1, \ldots, n.
\end{equation}

\vspace{2mm}

\subsection{Fundamental structure conditions}~

\vspace{2mm}

Let $K \subset \mathbb{R}^n$ denote an open symmetric convex cone and assume $\Gamma_n \subset K$.
As in \cite{GS10}, we assume the curvature function $f$ satisfies the following fundamental structure conditions.
\begin{equation} \label{eq2-16}
f_i (\lambda) \equiv \frac{\partial f}{\partial \lambda_i} > 0 \quad \text{ in } K, \ \ 1 \leq i \leq n,
\end{equation}
\begin{equation} \label{eq2-17}
f \text{ is a concave function in } K,
\end{equation}
\begin{equation} \label{eq2-18}
f > 0 \quad \text{ in } K, \quad f = 0 \quad \text{ on }  \partial K,
\end{equation}
\begin{equation} \label{eq2-20}
f \text{ is homogeneous of degree one},
\end{equation}
\begin{equation} \label{eq2-19}
f(1, \cdots, 1) = 1,
\end{equation}
\begin{equation} \label{eq2-21}
\lim\limits_{R \rightarrow \infty} f(\lambda_1, \ldots, \lambda_{n - 1}, \lambda_n + R) \geq 1 + \epsilon_0 \quad \text{uniformly in } B_{\delta_0}({\bf 1})
\end{equation}
for some fixed $\epsilon_0 > 0$ and $\delta_0 > 0$, where $B_{\delta_0}({\bf 1})$ is the ball of radius $\delta_0$ centered at ${\bf 1} = (1, \cdots, 1) \in \mathbb{R}^n$.

Denoting the vector space of $n \times n$ symmetric matrices by $\mathcal{S}$, and setting
\[ \mathcal{S}_K = \{ A \in \mathcal{S} | \lambda(A) \in K \}, \]
where $\lambda(A) = (\lambda_1, \cdots, \lambda_n)$ denotes the eigenvalues of $A$,
we define a function $F$ by
\[ F ( A ) = f \big(  \lambda ( A ) \big), \quad A \in \mathcal{S}_K. \]
Denote 
\[ F^{ij}(A) = \frac{\partial F}{\partial a_{ij}} (A), \quad F^{ij, kl} (A) = \frac{\partial^2 F}{\partial a_{ij} \partial a_{kl}} (A) . \]
Next, we define a function $G$ by
\[ G(D^2 u, D u, u) = F ( A ), \]
where $A$ is given by \eqref{eq2-13}.
Now we can see that equation \eqref{eq1-1} can be rewritten as
\begin{equation} \label{eq2-14}
G(D^2 u, D u, u) = \sigma.
\end{equation}

\vspace{2mm}

\subsection{Calculations on hypersurface $\Sigma$}~

\vspace{2mm}

Let $\tilde{g}$ and $\tilde{\nabla}$ denote the metric and Levi-Civita connection induced from $\mathbb{R}^{n + 1}$, and $g$ and $\nabla$ denote the metric and Levi-Civita connection on $\Sigma$ induced from $\mathbb{H}^{n + 1}$. The Christoffel symbols with respect to $\nabla$ and $\tilde{\nabla}$ have the following relation
\[ \Gamma_{ij}^k = \tilde{\Gamma}_{ij}^k - \frac{1}{u} (u_i \delta_{kj} + u_j \delta_{ik} - \tilde{g}^{kl} u_l \tilde{g}_{ij}). \]
Thus, for any $C^2$ function $v$ defined on $\Sigma$, we know that
\begin{equation} \label{eq2-1}
\nabla_{ij} v = (v_i)_j - \Gamma_{ij}^k v_k = \tilde{\nabla}_{ij} v + \frac{1}{u}( u_i v_j + u_j v_i - \tilde{g}^{kl} u_l v_k \tilde{g}_{ij} ).
\end{equation}
We shall need the following lemmas.
\begin{lemma}  \label{Lemma1}
For $\Sigma \subset \mathbb{R}^{n+1}$, the following identities hold.
\begin{equation} \label{eq2-2}
\tilde{g}^{kl} u_k u_l  = |\tilde\nabla u|^2 = 1 - (\nu^{n + 1})^2,
\end{equation}
\begin{equation}  \label{eq2-3}
\tilde{\nabla}_{ij} u = \tilde{h}_{ij} \nu^{n+1} \quad \mbox{and} \quad \tilde{\nabla}_{ij} x_{k} = \tilde{h}_{ij} \nu^{k}, \quad k = 1, \ldots, n,
\end{equation}
\begin{equation}  \label{eq2-4}
(\nu^{n+1})_i = - \tilde{h}_{ij} \tilde{g}^{j k} u_k,
\end{equation}
\begin{equation} \label{eq2-5}
\tilde{\nabla}_{ij} \nu^{n+1} = - \tilde{g}^{kl} ( \nu^{n+1} \tilde{h}_{il} \tilde{h}_{kj} + u_l \tilde{\nabla}_k \tilde{h}_{ij} ) .
\end{equation}
\end{lemma}

\begin{lemma}  \label{Lemma6}
Let $\Sigma$ be a smooth hypersurface in $\mathbb{H}^{n+1}$ which satisfies equation \eqref{eq1-1}. Then in a local orthonormal frame on $\Sigma$, we have the identity
\begin{equation}  \label{eq4-1}
\begin{aligned}
F^{ij} \nabla_{ij} \nu^{n+1}
= &  \Big( 1 + (\nu^{n+1})^2 \Big) \sigma - \nu^{n+1} \Big( \sum f_i +  \sum f_i \kappa_i^2 \Big) \\ & + \frac{2}{u^2} F^{ij} u_i u_j \big( \nu^{n + 1} - \kappa_j \big).
\end{aligned}
\end{equation}
\end{lemma}

\vspace{4mm}

\section{Uniform global curvature estimate}

\vspace{4mm}

For curvature function \eqref{eq1-6}, direct calculation shows that
\begin{equation} \label{eq5-1}
f_i = \frac{\partial f}{\partial \kappa_i} = \frac{f}{C_n^k} \sum\limits_{\substack{ i_1 < \cdots < i_k \\ i \in \{ i_1, \cdots, i_k \} } } \frac{1}{\kappa_{i_1} + \cdots + \kappa_{i_k}}.
\end{equation}

In this section, we begin uniform global curvature estimate. Since $\Gamma$ is assumed to be mean convex, by Guan and Spruck \cite{GS10} we know that
\[ \nu^{n + 1} \geq \sigma   \quad \text{on} \quad \Sigma . \]
We consider the test function
\begin{equation*}
M_0 = \sup\limits_{X \in \Sigma} \frac{H}{({\nu}^{n+1})^{\beta}},
\end{equation*}
where $H = \sum \kappa_i$ is the mean curvature of $\Sigma$ and $\beta > 1$ is a constant to be determined later.

Assume that $M_0 > 0$ is attained at an interior point $X_0 \in \Sigma$.
Let $\tau_1, \ldots, \tau_n$ be a smooth local orthonormal frame field about
$X_0$ such that $h_{ij}(X_0) = \kappa_i \delta_{ij}$, where
$\kappa_1, \ldots, \kappa_n$ are the hyperbolic principal curvatures of
$\Sigma$ at $X_0$. We assume $\kappa_1 \geq \cdots \geq \kappa_n$.
Since
\[ \ln H - \beta \ln ( {\nu}^{n+1} ) \]
also achieves its maximum at $X_0$, therefore at $X_0$, we have
\begin{equation} \label{eq4-2}
\frac{\sum_j h_{jji}}{H} - \frac{\beta \nabla_i \nu^{n + 1}}{\nu^{n + 1}} = 0,
\end{equation}

\begin{equation} \label{eq4-3}
\frac{\sum_j h_{jjii}}{H} - \frac{\beta \nabla_{ii} \nu^{n + 1}}{\nu^{n + 1}} - \frac{\beta (\beta - 1)( \nabla_i \nu^{n + 1} )^2}{(\nu^{n + 1})^2} \leq 0.
\end{equation}

Now we differentiate equation \eqref{eq1-1} twice to obtain
\begin{equation} \label{eq4-23}
\sum_i F^{ii} h_{iij} = 0,
\end{equation}
and
\begin{equation}  \label{eq4-4}
\sum_i F^{ii} h_{iijj} + \sum_{p q r s} F^{pq, rs} h_{pqj} h_{rsj} = 0 .
\end{equation}
Besides, Gauss equation implies
\begin{equation} \label{eq2G-4}
h_{iijj} = h_{jjii} + ( \kappa_i \kappa_j - 1 )( \kappa_i -
\kappa_j ).
\end{equation}

Denoting
\[ A = \kappa_1^2 + \cdots + \kappa_n^2, \]
and combining \eqref{eq4-3}, \eqref{eq4-4}, \eqref{eq2G-4} and \eqref{eq4-1} yields
\begin{equation} \label{eq4-5}
\begin{aligned}
& \bigg( A - \frac{\beta \big( 1 + (\nu^{n+1})^2 \big)}{\nu^{n + 1}} H + n \bigg) \sigma + (\beta - 1) H \bigg( \sum f_i + \sum f_i \kappa_i^2 \bigg) \\
& - \sum_{pqrsj} F^{pq, rs} h_{pqj} h_{rsj}
+ \frac{2 \beta}{\nu^{n + 1}} H \sum f_i \frac{u_i^2}{u^2} \big( \kappa_i - \nu^{n+1} \big) \\
& - \frac{\beta (\beta - 1)}{(\nu^{n + 1})^2} H \sum f_i ( \nabla_i \nu^{n + 1} )^2 \leq 0.
\end{aligned}
\end{equation}

By \eqref{eq2-4} and \eqref{eq2-11}, we have
\begin{equation} \label{eq4-24}
(\nu^{n + 1})_i = \frac{u_i}{u} (\nu^{n + 1} - \kappa_i).
\end{equation}
Taking \eqref{eq4-24} into \eqref{eq4-2} yields
\begin{equation} \label{eq4-6}
\sum_j h_{jji} = \frac{\beta}{\nu^{n + 1}} \frac{u_i}{u} (\nu^{n+1} - \kappa_i) H.
\end{equation}

We note that
\begin{equation} \label{eq4-7}
\begin{aligned}
& - \sum_{pqrs} F^{pq, rs} h_{pqj} h_{rsj}  =  \sum\limits_{p \neq q} \frac{f_p - f_q}{\kappa_q - \kappa_p} h_{pqj}^2 - \sum\limits_{pq} \frac{\partial^2 f}{\partial \kappa_p \partial \kappa_q} h_{ppj} h_{qqj}.
\end{aligned}
\end{equation}
One may refer to \cite{Ball1984, GRW2015} to find this formula.

Taking \eqref{eq4-24} and \eqref{eq4-7} into \eqref{eq4-5} yields 

\begin{equation} \label{eq4-8}
\begin{aligned}
& \bigg( A - \frac{\beta \big( 1 + (\nu^{n+1})^2 \big)}{\nu^{n + 1}} H + n \bigg) \sigma + (\beta - 1) H \bigg( \sum f_i + \sum f_i \kappa_i^2 \bigg) \\
&  + \sum_j \sum\limits_{p \neq q} \frac{f_p - f_q}{\kappa_q - \kappa_p} h_{pqj}^2 - \sum\limits_{pqj} \frac{\partial^2 f}{\partial \kappa_p \partial \kappa_q} h_{ppj} h_{qqj}
\\
& + \frac{2 \beta}{\nu^{n + 1}} H \sum f_i \frac{u_i^2}{u^2} \big( \kappa_i - \nu^{n+1} \big) - \frac{\beta (\beta - 1)}{(\nu^{n + 1})^2} H \sum f_i  \frac{u_i^2}{u^2} (\nu^{n + 1} - \kappa_i)^2 \leq 0.
\end{aligned}
\end{equation}
Then we take 
\begin{equation} \label{eq4-10}
\sum_j \sum\limits_{p \neq q} \frac{f_p - f_q}{\kappa_q - \kappa_p} h_{pqj}^2
\geq 2 \sum\limits_{i \neq j} \frac{f_i - f_j}{\kappa_j - \kappa_i} h_{ijj}^2
\end{equation}
and
\begin{equation}   \label{eq4-11}
 - \sum\limits_{pqj} \frac{\partial^2 f}{\partial \kappa_p \partial \kappa_q} h_{ppj} h_{qqj}
= - \sum\limits_{i j l} f_{jl}  h_{jji} h_{lli}
\end{equation}
into \eqref{eq4-8} to obtain 
\begin{equation} \label{eq4-9}
\begin{aligned}
& \bigg( A - \frac{\beta \big( 1 + (\nu^{n+1})^2 \big)}{\nu^{n + 1}} H + n \bigg) \sigma + ( \beta - 1 ) H \bigg( \sum f_i + \sum f_i \kappa_i^2 \bigg) \\
&  + 2 \sum\limits_{i \neq j} \frac{f_i - f_j}{\kappa_j - \kappa_i} h_{ijj}^2  - \sum\limits_{i j l} f_{jl}  h_{jji} h_{lli}
\\
& + \frac{2 \beta}{\nu^{n + 1}} H \sum f_i \frac{u_i^2}{u^2} \big( \kappa_i - \nu^{n+1} \big) - \frac{\beta (\beta - 1)}{(\nu^{n + 1})^2} H \sum f_i  \frac{u_i^2}{u^2} (\nu^{n + 1} - \kappa_i)^2 \leq 0.
\end{aligned}
\end{equation}

\vspace{4mm}

\section{Computing concavity by Lagrange multiplier method}

\vspace{4mm}

Hereafter, we confine ourselves to the case when $k = n - 1$ in \eqref{eq1-6}.
In this case,
\begin{equation*}
f_j = \frac{\sigma}{n} \sum\limits_{p \neq j} \frac{1}{H - \kappa_p},
\end{equation*}

\begin{equation*}
\begin{aligned}
f_{jl} = & \left\{\begin{aligned}
& - \frac{\sigma}{n} \sum\limits_{p \neq j} \frac{1}{(H - \kappa_p)^2} + \frac{f_j^2}{\sigma}, \quad j = l, \\
& - \frac{\sigma}{n} \sum\limits_{p \notin \{ j, l \}} \frac{1}{(H - \kappa_p)^2} + \frac{f_j f_l}{\sigma}, \quad j \neq l
\end{aligned}
\right. \\
= &  \left\{\begin{aligned}
& \Big( - \frac{\sigma}{n} + \frac{\sigma}{n^2} \Big) \sum\limits_{p \neq j} \frac{1}{(H - \kappa_p)^2} + \frac{2 \sigma}{n^2} \sum\limits_{\substack{p < q \\ p, q \neq j }} \frac{1}{(H - \kappa_p)(H - \kappa_q)} , \quad j = l, \\
&  \Big( - \frac{\sigma}{n} + \frac{\sigma}{n^2} \Big) \sum\limits_{p \notin \{ j, l \}} \frac{1}{(H - \kappa_p)^2} + \frac{\sigma}{n^2} \sum\limits_{\substack{p < q \\ p, q \notin \{j, l\} }} \frac{1}{(H - \kappa_p)(H - \kappa_q)} \\
& + \frac{\sigma}{n^2} \sum\limits_{p < q} \frac{1}{(H - \kappa_p)(H - \kappa_q)}, \quad j \neq l.
\end{aligned}
\right.
\end{aligned}
\end{equation*}

When $j \neq i$,
\begin{equation*}
 \frac{f_i - f_j}{\kappa_j - \kappa_i} = \frac{\sigma}{n} \frac{1}{(H - \kappa_i)(H - \kappa_j)}.
\end{equation*}

We note that
\begin{equation*}
f_j f_l = \left\{
\begin{aligned}
& \frac{\sigma^2}{n^2} \sum\limits_{p \neq j} \frac{1}{(H - \kappa_p)^2} + \frac{2 \sigma^2}{n^2} \sum\limits_{\substack{p < q \\ p, q \neq j }} \frac{1}{(H - \kappa_p)(H - \kappa_q)} , \quad j = l, \\
& \frac{\sigma^2}{n^2} \sum\limits_{p \notin \{ j, l \}} \frac{1}{(H - \kappa_p)^2} + \frac{\sigma^2}{n^2} \sum\limits_{\substack{p < q \\ p, q \notin \{j, l\} }} \frac{1}{(H - \kappa_p)(H - \kappa_q)} \\
& + \frac{\sigma^2}{n^2} \sum\limits_{p < q} \frac{1}{(H - \kappa_p)(H - \kappa_q)}, \quad j \neq l.
\end{aligned} \right.
\end{equation*}

For each $i = 1, \ldots, n$, we shall minimize the quantity
\begin{equation} \label{eq4-12}
Q_i := 2 \sum\limits_{j \neq i} \frac{f_i - f_j}{\kappa_j - \kappa_i} h_{ijj}^2  - \sum\limits_{j l} f_{jl}  h_{jji} h_{lli} .
\end{equation}

By \eqref{eq4-23} and \eqref{eq4-6} we can see that \eqref{eq4-12} becomes
\begin{equation*}
\begin{aligned}
Q_i = & 2 \sum\limits_{j \neq i} \frac{f_i - f_j}{\kappa_j - \kappa_i} h_{ijj}^2  - \sum\limits_{j l} \Big( f_{jl} + \frac{n - 1}{\sigma} f_j f_l \Big)  h_{jji} h_{lli} \\
= &  2 \sum\limits_{j \neq i} \frac{f_i - f_j}{\kappa_j - \kappa_i} h_{ijj}^2  - \sum\limits_{j l} \Big( f_{jl} + \frac{n - 1}{\sigma} f_j f_l - \frac{\sigma}{n} \sum\limits_{p < q} \frac{1}{(H - \kappa_p)(H - \kappa_q)} \Big)  h_{jji} h_{lli} \\
& - \frac{\beta^2}{(\nu^{n + 1})^2} \frac{u_i^2}{u^2} (\nu^{n+1} - \kappa_i)^2 H^2 \frac{\sigma}{n} \sum\limits_{p < q} \frac{1}{(H - \kappa_p)(H - \kappa_q)} \\
= &  \frac{2 \sigma}{n} \frac{1}{H - \kappa_i} \sum\limits_{j \neq i} \frac{1}{H - \kappa_j} h_{ijj}^2  - \frac{\sigma}{n} \sum\limits_{j \neq l} \sum\limits_{\substack{p < q \\ p, q \notin \{j, l\} }} \frac{1}{(H - \kappa_p)(H - \kappa_q)} h_{jji} h_{lli} \\
& - \frac{\sigma}{n} \sum\limits_j \bigg( 2 \sum\limits_{\substack{p < q \\ p, q \neq j }} \frac{1}{(H - \kappa_p)(H - \kappa_q)} -  \sum\limits_{p < q} \frac{1}{(H - \kappa_p)(H - \kappa_q)} \bigg) h_{jji}^2 \\
& - \frac{\beta^2}{(\nu^{n + 1})^2} \frac{u_i^2}{u^2} (\nu^{n+1} - \kappa_i)^2 H^2 \frac{\sigma}{n} \sum\limits_{p < q} \frac{1}{(H - \kappa_p)(H - \kappa_q)} .
\end{aligned}
\end{equation*}

Now we shall minimize $Q_i$ under condition \eqref{eq4-6} and \eqref{eq4-23}. It suffices to minimize the function $L_i$ in terms of $t_1$, \ldots, $t_n$, $\mu_1$, $\mu_2$,

\begin{equation*}
\begin{aligned}
L_i = &  \frac{2 \sigma}{n} \frac{1}{H - \kappa_i} \sum\limits_{j \neq i} \frac{1}{H - \kappa_j} t_j^2  - \frac{\sigma}{n} \sum\limits_{j \neq l} \sum\limits_{\substack{p < q \\ p, q \notin \{j, l\} }} \frac{1}{(H - \kappa_p)(H - \kappa_q)} t_j t_l \\
& - \frac{\sigma}{n} \sum\limits_j \bigg( 2 \sum\limits_{\substack{p < q \\ p, q \neq j }} \frac{1}{(H - \kappa_p)(H - \kappa_q)} -  \sum\limits_{p < q} \frac{1}{(H - \kappa_p)(H - \kappa_q)} \bigg) t_j^2 \\
& - \frac{\beta^2}{(\nu^{n + 1})^2} \frac{u_i^2}{u^2} (\nu^{n+1} - \kappa_i)^2 H^2 \frac{\sigma}{n} \sum\limits_{p < q} \frac{1}{(H - \kappa_p)(H - \kappa_q)} \\
& - \mu_1  \sum_j f_j t_j - \mu_2 \bigg( \sum_j t_j - \frac{\beta}{\nu^{n + 1}} \frac{u_i}{u} (\nu^{n+1} - \kappa_i) H \bigg).
\end{aligned}
\end{equation*}

To find the critical points of $L_i$, we solve the following linear equations.

\begin{equation} \label{eq3-8}
\begin{aligned}
\frac{\partial L_i}{\partial t_j} = &  \frac{4 \sigma}{n} \frac{1}{(H - \kappa_i)(H - \kappa_j)} t_j  - \frac{2 \sigma}{n} \sum\limits_{l \neq j} \sum\limits_{\substack{p < q \\ p, q \notin \{j, l\} }} \frac{1}{(H - \kappa_p)(H - \kappa_q)} t_l \\
& - \frac{2 \sigma}{n}  \bigg( 2 \sum\limits_{\substack{p < q \\ p, q \neq j }} \frac{1}{(H - \kappa_p)(H - \kappa_q)} -  \sum\limits_{p < q} \frac{1}{(H - \kappa_p)(H - \kappa_q)} \bigg) t_j \\
& - \mu_1  f_j - \mu_2 = 0, \qquad j \neq i, \qquad \qquad  & (j) \\
\frac{\partial L_i}{\partial t_i} = & - \frac{2 \sigma}{n} \sum\limits_{l \neq i} \sum\limits_{\substack{p < q \\ p, q \notin \{i, l\} }} \frac{1}{(H - \kappa_p)(H - \kappa_q)} t_l \\
& - \frac{2 \sigma}{n}  \bigg( 2 \sum\limits_{\substack{p < q \\ p, q \neq i }} \frac{1}{(H - \kappa_p)(H - \kappa_q)} -  \sum\limits_{p < q} \frac{1}{(H - \kappa_p)(H - \kappa_q)} \bigg) t_i \\
& - \mu_1  f_i - \mu_2 = 0,  \qquad \qquad & (i) \\
\sum_j f_j t_j = & 0,  \qquad \qquad & (d) \\
\sum_j t_j = & \frac{\beta}{\nu^{n + 1}} \frac{u_i}{u} (\nu^{n+1} - \kappa_i) H.  \qquad \qquad & (c)
\end{aligned}
\end{equation}

For $j \neq i$, $(j) - (i)$ yields,
\begin{equation} \label{eq3-7}
\begin{aligned}
& A_j^i t_i + A_j^j  t_j + \sum\limits_{l \neq i, j} A_j^l t_l =  \mu_1 ( f_j - f_i ),  \qquad \qquad & (j*)
\end{aligned}
\end{equation}
where
\begin{equation*}
\begin{aligned}
A_j^i = & \frac{4 \sigma}{n} \sum\limits_{\substack{p < q \\ p, q \neq i }} \frac{1}{(H - \kappa_p)(H - \kappa_q)} - \frac{2 \sigma}{n}  \sum\limits_{p < q} \frac{1}{(H - \kappa_p)(H - \kappa_q)} \\
& - \frac{2 \sigma}{n} \sum\limits_{\substack{p < q \\ p, q \notin \{j, i\} }} \frac{1}{(H - \kappa_p)(H - \kappa_q)},
\end{aligned}
\end{equation*}
\begin{equation*}
\begin{aligned}
A_j^j = &  \frac{4 \sigma}{n} \frac{1}{(H - \kappa_i)(H - \kappa_j)} - \frac{4 \sigma}{n} \sum\limits_{\substack{p < q \\ p, q \neq j }} \frac{1}{(H - \kappa_p)(H - \kappa_q)} \\
& + \frac{2 \sigma}{n} \sum\limits_{p < q} \frac{1}{(H - \kappa_p)(H - \kappa_q)} + \frac{2 \sigma}{n} \sum\limits_{\substack{p < q \\ p, q \notin \{i, j\} }} \frac{1}{(H - \kappa_p)(H - \kappa_q)},
\end{aligned}
\end{equation*}
\begin{equation*}
\begin{aligned}
A_j^l = \frac{2 \sigma}{n} \bigg( \sum\limits_{\substack{p < q \\ p, q \notin \{i, l\} }} \frac{1}{(H - \kappa_p)(H - \kappa_q)} - \sum\limits_{\substack{p < q \\ p, q \notin \{j, l\} }} \frac{1}{(H - \kappa_p)(H - \kappa_q)} \bigg), \qquad l \neq i, j.
\end{aligned}
\end{equation*}
Now we simplify $A_j^i$ and $A_j^j$.
We notice that
\begin{equation} \label{eq3-1}
\begin{aligned}
& \sum\limits_{\substack{p < q \\ p, q \neq i }} \frac{1}{(H - \kappa_p)(H - \kappa_q)} \\
= &  \sum\limits_{\substack{p < q \\ p, q \notin \{ i, j \} }} \frac{1}{(H - \kappa_p)(H - \kappa_q)} +   \sum\limits_{\substack{j < q \\  q \neq  i  }} \frac{1}{(H - \kappa_j)(H - \kappa_q)} +   \sum\limits_{\substack{p < j \\ p \neq i }} \frac{1}{(H - \kappa_p)(H - \kappa_j)} \\
= &  \sum\limits_{\substack{p < q \\ p, q \notin \{ i, j \} }} \frac{1}{(H - \kappa_p)(H - \kappa_q)} +   \frac{1}{H - \kappa_j} \sum\limits_{ p \neq i, j } \frac{1}{H - \kappa_p}.
\end{aligned}
\end{equation}
It follows that
\begin{equation} \label{eq3-2}
\begin{aligned}
& \sum\limits_{p < q} \frac{1}{(H - \kappa_p)(H - \kappa_q)} \\
= &  \sum\limits_{\substack{p < q \\ p, q \neq i} } \frac{1}{(H - \kappa_p)(H - \kappa_q)} + \sum\limits_{i < q} \frac{1}{(H - \kappa_i)(H - \kappa_q)} + \sum\limits_{p < i} \frac{1}{(H - \kappa_p)(H - \kappa_i)} \\
= &  \sum\limits_{\substack{p < q \\ p, q \neq i} } \frac{1}{(H - \kappa_p)(H - \kappa_q)} + \frac{1}{H - \kappa_i} \sum\limits_{p \neq i} \frac{1}{H - \kappa_p} \\
= &   \sum\limits_{\substack{p < q \\ p, q \notin \{ i, j \} }} \frac{1}{(H - \kappa_p)(H - \kappa_q)} +   \frac{1}{H - \kappa_j} \sum\limits_{ p \neq i, j } \frac{1}{H - \kappa_p} \\
& + \frac{1}{H - \kappa_i} \sum\limits_{p \neq i, j} \frac{1}{H - \kappa_p} + \frac{1}{H - \kappa_i} \frac{1}{H - \kappa_j}.
\end{aligned}
\end{equation}
By \eqref{eq3-1} and \eqref{eq3-2},
\begin{equation} \label{eq3-3}
\begin{aligned}
A_j^i = & \frac{4 \sigma}{n} \sum\limits_{\substack{p < q \\ p, q \notin \{ i, j \} }} \frac{1}{(H - \kappa_p)(H - \kappa_q)} +  \frac{4 \sigma}{n} \frac{1}{H - \kappa_j} \sum\limits_{ p \neq i, j } \frac{1}{H - \kappa_p} \\
 & - \frac{2 \sigma}{n} \sum\limits_{\substack{p < q \\ p, q \notin \{ i, j \} }} \frac{1}{(H - \kappa_p)(H - \kappa_q)} - \frac{2 \sigma}{n} \frac{1}{H - \kappa_j} \sum\limits_{ p \neq i, j } \frac{1}{H - \kappa_p} \\
& - \frac{2 \sigma}{n} \frac{1}{H - \kappa_i} \sum\limits_{p \neq i, j} \frac{1}{H - \kappa_p} - \frac{2 \sigma}{n} \frac{1}{H - \kappa_i} \frac{1}{H - \kappa_j} \\
& - \frac{2 \sigma}{n} \sum\limits_{\substack{p < q \\ p, q \notin \{j, i\} }} \frac{1}{(H - \kappa_p)(H - \kappa_q)} \\
= &  \frac{2 \sigma}{n} \frac{\kappa_j - \kappa_i}{(H - \kappa_j)(H - \kappa_i)} \sum\limits_{ p \neq i, j } \frac{1}{H - \kappa_p}  - \frac{2 \sigma}{n} \frac{1}{(H - \kappa_i)(H - \kappa_j)}.
\end{aligned}
\end{equation}
Similarly,
\begin{equation} \label{eq3-5}
\begin{aligned}
& \sum\limits_{\substack{p < q \\ p, q \neq j }} \frac{1}{(H - \kappa_p)(H - \kappa_q)} \\
= &  \sum\limits_{\substack{p < q \\ p, q \notin \{ j, i \} }} \frac{1}{(H - \kappa_p)(H - \kappa_q)} +   \sum\limits_{\substack{i < q \\  q \neq  j  }} \frac{1}{(H - \kappa_i)(H - \kappa_q)} +   \sum\limits_{\substack{p < i \\ p \neq j }} \frac{1}{(H - \kappa_p)(H - \kappa_i)} \\
= &  \sum\limits_{\substack{p < q \\ p, q \notin \{ i, j \} }} \frac{1}{(H - \kappa_p)(H - \kappa_q)} +   \frac{1}{H - \kappa_i} \sum\limits_{ p \neq i, j } \frac{1}{H - \kappa_p}.
\end{aligned}
\end{equation}
\begin{equation} \label{eq3-4}
\begin{aligned}
A_j^j = &  \frac{4 \sigma}{n} \frac{1}{(H - \kappa_i)(H - \kappa_j)} \\
& - \frac{4 \sigma}{n} \sum\limits_{\substack{p < q \\ p, q \notin \{ i, j \} }} \frac{1}{(H - \kappa_p)(H - \kappa_q)} - \frac{4 \sigma}{n} \frac{1}{H - \kappa_i} \sum\limits_{ p \neq i, j } \frac{1}{H - \kappa_p} \\
& + \frac{2 \sigma}{n} \sum\limits_{\substack{p < q \\ p, q \notin \{ i, j \} }} \frac{1}{(H - \kappa_p)(H - \kappa_q)} + \frac{2 \sigma}{n}  \frac{1}{H - \kappa_j} \sum\limits_{ p \neq i, j } \frac{1}{H - \kappa_p} \\
& + \frac{2 \sigma}{n} \frac{1}{H - \kappa_i} \sum\limits_{p \neq i, j} \frac{1}{H - \kappa_p} + \frac{2 \sigma}{n} \frac{1}{H - \kappa_i} \frac{1}{H - \kappa_j} \\
& + \frac{2 \sigma}{n} \sum\limits_{\substack{p < q \\ p, q \notin \{i, j\} }} \frac{1}{(H - \kappa_p)(H - \kappa_q)} \\
= &  \frac{6 \sigma}{n} \frac{1}{(H - \kappa_i)(H - \kappa_j)} + \frac{2 \sigma}{n} \frac{\kappa_j - \kappa_i}{(H - \kappa_i)(H - \kappa_j)} \sum\limits_{ p \neq i, j } \frac{1}{H - \kappa_p}.
\end{aligned}
\end{equation}

$A_j^l$ can be simplified as
\begin{equation} \label{eq3-6}
\begin{aligned}
A_j^l = &  \frac{2 \sigma}{n} \sum\limits_{\substack{p < q \\ p, q \notin \{i, l, j\} }} \frac{1}{(H - \kappa_p)(H - \kappa_q)} + \frac{2 \sigma}{n} \sum\limits_{\substack{j < q \\ q \neq i, l }} \frac{1}{(H - \kappa_j)(H - \kappa_q)} \\
& + \frac{2 \sigma}{n} \sum\limits_{\substack{p < j \\ p \neq i, l }} \frac{1}{(H - \kappa_p)(H - \kappa_j)}  \\
& - \frac{2 \sigma}{n} \sum\limits_{\substack{p < q \\ p, q \notin \{j, l, i\} }} \frac{1}{(H - \kappa_p)(H - \kappa_q)}  - \frac{2 \sigma}{n} \sum\limits_{\substack{i < q \\  q \notin \{j, l\} }} \frac{1}{(H - \kappa_i)(H - \kappa_q)} \\
& - \frac{2 \sigma}{n} \sum\limits_{\substack{p < i \\ p \neq j, l }} \frac{1}{(H - \kappa_p)(H - \kappa_i)}   \\
= & \frac{2 \sigma}{n} \frac{1}{H - \kappa_j} \sum\limits_{ p \neq i, l, j } \frac{1}{H - \kappa_p} - \frac{2 \sigma}{n} \frac{1}{H - \kappa_i} \sum\limits_{ p \neq j, l, i } \frac{1}{H - \kappa_p} \\
= &  \frac{2 \sigma}{n} \frac{\kappa_j - \kappa_i}{(H - \kappa_i)(H - \kappa_j)} \sum\limits_{ p \neq i, l, j } \frac{1}{H - \kappa_p},    & \qquad l \neq i, j.
\end{aligned}
\end{equation}
By \eqref{eq3-3}, \eqref{eq3-4} and \eqref{eq3-6}, we know that $(j*)$ in \eqref{eq3-7} becomes
\begin{equation*}
\begin{aligned}
& \bigg( \frac{2 \sigma}{n} \frac{\kappa_j - \kappa_i}{(H - \kappa_j)(H - \kappa_i)} \sum\limits_{ p \neq i, j } \frac{1}{H - \kappa_p}  - \frac{2 \sigma}{n} \frac{1}{(H - \kappa_i)(H - \kappa_j)} \bigg) t_i \\
& + \bigg( \frac{6 \sigma}{n} \frac{1}{(H - \kappa_i)(H - \kappa_j)} + \frac{2 \sigma}{n} \frac{\kappa_j - \kappa_i}{(H - \kappa_i)(H - \kappa_j)} \sum\limits_{ p \neq i, j } \frac{1}{H - \kappa_p} \bigg) t_j \\
& + \sum\limits_{l \neq i, j}   \frac{2 \sigma}{n} \frac{\kappa_j - \kappa_i}{(H - \kappa_i)(H - \kappa_j)} \sum\limits_{ p \neq i, l, j } \frac{1}{H - \kappa_p} t_l =  \mu_1 \frac{\sigma}{n} \frac{\kappa_i - \kappa_j}{(H - \kappa_i)(H - \kappa_j)},
\end{aligned}
\end{equation*}
which is equivalent to
\begin{equation} \label{eq3-9}
\begin{aligned}
& \bigg( (\kappa_j - \kappa_i) \sum\limits_{ p \neq i, j } \frac{1}{H - \kappa_p}  - 1  \bigg) t_i  + \bigg( 3 + (\kappa_j - \kappa_i) \sum\limits_{ p \neq i, j } \frac{1}{H - \kappa_p} \bigg) t_j \\
& + \sum\limits_{l \neq i, j} (\kappa_j - \kappa_i) \sum\limits_{ p \neq i, l, j } \frac{1}{H - \kappa_p} t_l = \frac{\mu_1}{2} (\kappa_i - \kappa_j), \qquad j \neq i.  \qquad \qquad & (j**)
\end{aligned}
\end{equation}
Adding all $(j**)$ for $j \neq i$ yields
\begin{equation} \label{eq3-11}
\begin{aligned}
& \bigg( \sum\limits_{j \neq i} (\kappa_j - \kappa_i) \sum\limits_{ p \neq i, j } \frac{1}{H - \kappa_p}  - (n - 1) \bigg) t_i  + \sum\limits_{j \neq i} \bigg( 3 + (\kappa_j - \kappa_i) \sum\limits_{ p \neq i, j } \frac{1}{H - \kappa_p} \bigg) t_j \\
& + \sum\limits_{j \neq i} \sum\limits_{l \neq i, j} (\kappa_j - \kappa_i) \sum\limits_{ p \neq i, l, j } \frac{1}{H - \kappa_p} t_l = \frac{\mu_1}{2} \sum\limits_{j \neq i} (\kappa_i - \kappa_j).
\end{aligned}
\end{equation}
We note that
\begin{equation} \label{eq3-10}
\begin{aligned}
& \sum\limits_{j \neq i} \sum\limits_{l \neq i, j} (\kappa_j - \kappa_i) \sum\limits_{ p \neq i, l, j } \frac{1}{H - \kappa_p} t_l \\
= & \sum\limits_{l \neq i} \sum\limits_{p \neq i, l} \frac{1}{H - \kappa_p} \sum\limits_{ j \neq i, l, p } (\kappa_j - \kappa_i) t_l \\
= & \sum\limits_{l \neq i} \sum\limits_{p \neq i, l} \frac{1}{H - \kappa_p} \Big( H - (n - 2) \kappa_i - \kappa_p - \kappa_l \Big) t_l \\
= & \sum\limits_{l \neq i} \sum\limits_{p \neq i, l}  \bigg( 1 - \frac{(n - 2) \kappa_i + \kappa_l}{H - \kappa_p}  \bigg) t_l \\
= & \sum\limits_{l \neq i}  \bigg( n - 2 - \Big( (n - 2) \kappa_i + \kappa_l \Big) \sum\limits_{p \neq i, l} \frac{1}{H - \kappa_p}  \bigg) t_l \\
= & \sum\limits_{j \neq i}  \bigg( n - 2 - \Big( (n - 2) \kappa_i + \kappa_j \Big) \sum\limits_{p \neq i, j} \frac{1}{H - \kappa_p}  \bigg) t_j.
\end{aligned}
\end{equation}
Besides,
\begin{equation} \label{eq3-12}
\begin{aligned}
& \sum\limits_{j \neq i} (\kappa_j - \kappa_i) \sum\limits_{ p \neq i, j } \frac{1}{H - \kappa_p}  - (n - 1) \\
= & \sum\limits_{p \neq i} \frac{1}{H - \kappa_p} \sum\limits_{ j \neq i, p } (\kappa_j - \kappa_i) - (n - 1) \\
= & \sum\limits_{p \neq i} \frac{1}{H - \kappa_p} \Big( H - (n - 1) \kappa_i - \kappa_p \Big) - (n - 1) \\
= & \sum\limits_{p \neq i} \bigg( 1 - \frac{(n - 1) \kappa_i}{H - \kappa_p} \bigg) - (n - 1) \\
= & - (n - 1) \kappa_i \sum\limits_{p \neq i} \frac{1}{H - \kappa_p} .
\end{aligned}
\end{equation}
Taking \eqref{eq3-10} and \eqref{eq3-12} into \eqref{eq3-11} we obtain
\begin{equation*}
\begin{aligned}
&  - (n - 1) \kappa_i \sum\limits_{p \neq i} \frac{1}{H - \kappa_p} t_i  + \sum\limits_{j \neq i} \bigg( n + 1  - (n - 1) \kappa_i \sum\limits_{ p \neq i, j } \frac{1}{H - \kappa_p} \bigg) t_j  = \frac{\mu_1}{2} (n \kappa_i - H),
\end{aligned}
\end{equation*}
which is equivalent to
\begin{equation}  \label{eq3-13}
\begin{aligned}
&  - (n - 1) \kappa_i f_i t_i  + (n + 1) \frac{\sigma}{n} \sum\limits_{j \neq i} t_j - (n - 1) \kappa_i \sum\limits_{j \neq i}   \bigg( f_j - \frac{\sigma}{n} \frac{1}{H - \kappa_i} \bigg) t_j = \frac{\mu_1}{2} \frac{\sigma}{n} (n \kappa_i - H).
\end{aligned}
\end{equation}

In view of equation $(d)$ and $(c)$ in \eqref{eq3-8}, \eqref{eq3-13} can be reduced to
\begin{equation*}
\begin{aligned}
& \bigg( (n + 1) + \frac{(n - 1) \kappa_i}{H - \kappa_i} \bigg) \bigg( \frac{\beta}{\nu^{n + 1}} \frac{u_i}{u} (\nu^{n+1} - \kappa_i) H - t_i \bigg) = \frac{\mu_1}{2} (n \kappa_i - H),
\end{aligned}
\end{equation*}
which implies that
\begin{equation}  \label{eq3-14}
\begin{aligned}
t_i = \frac{\beta}{\nu^{n + 1}} \frac{u_i}{u} (\nu^{n+1} - \kappa_i) H - \frac{\mu_1}{2} \frac{(n \kappa_i - H) (H - \kappa_i)}{(n + 1) H - 2 \kappa_i}.
\end{aligned}
\end{equation}

Now we look at $(j**)$ in \eqref{eq3-9}. For $j \neq i$, we note that $(j**)$ can be rewritten as
\begin{equation} \label{eq3-15}
\begin{aligned}
& \bigg( (\kappa_j - \kappa_i)  \Big( f_i - \frac{\sigma}{n} \frac{1}{H - \kappa_j} \Big) - \frac{\sigma}{n}  \bigg) t_i  + \bigg( \frac{3 \sigma}{n} + (\kappa_j - \kappa_i)  \Big( f_j - \frac{\sigma}{n} \frac{1}{H - \kappa_i} \Big) \bigg) t_j \\
& + \sum\limits_{l \neq i, j} (\kappa_j - \kappa_i) \Big( f_l - \frac{\sigma}{n} \frac{1}{H - \kappa_i} - \frac{\sigma}{n} \frac{1}{H - \kappa_j} \Big) t_l = \frac{\mu_1}{2} \frac{\sigma}{n} (\kappa_i - \kappa_j).
\end{aligned}
\end{equation}
In view of equation $(d)$ and $(c)$ in \eqref{eq3-8}, \eqref{eq3-15} can be further reduced to
\begin{equation*}
\begin{aligned}
& \bigg( - (\kappa_j - \kappa_i) \frac{1}{H - \kappa_j}  - 1 \bigg) t_i  + \bigg( 3 - (\kappa_j - \kappa_i) \frac{1}{H - \kappa_i}  \bigg) t_j \\
& - (\kappa_j - \kappa_i) \bigg( \frac{1}{H - \kappa_i} + \frac{1}{H - \kappa_j} \bigg) \bigg(  \frac{\beta}{\nu^{n + 1}} \frac{u_i}{u} (\nu^{n+1} - \kappa_i) H - t_i - t_j \bigg) \\
= & \frac{\mu_1}{2} (\kappa_i - \kappa_j),
\end{aligned}
\end{equation*}
which can be further simplified as
\begin{equation} \label{eq3-16}
\begin{aligned}
& \frac{\kappa_j - H}{H - \kappa_i} t_i  + \frac{3 H - 2 \kappa_j - \kappa_i}{H - \kappa_j} t_j \\
= & \frac{\mu_1}{2} (\kappa_i - \kappa_j) + (\kappa_j - \kappa_i) \frac{(2 H - \kappa_j - \kappa_i)}{(H - \kappa_i)(H - \kappa_j)} \frac{\beta}{\nu^{n + 1}} \frac{u_i}{u} (\nu^{n+1} - \kappa_i) H .
\end{aligned}
\end{equation}
Taking \eqref{eq3-14} into \eqref{eq3-16} we have
\begin{equation} \label{eq3-17}
\begin{aligned}
t_j
= & \frac{\mu_1}{2} \frac{H - \kappa_j}{3 H - 2 \kappa_j - \kappa_i} \frac{H \kappa_i - (n + 2) H \kappa_j - 2 \kappa_i^2 + (n + 2) \kappa_i \kappa_j + H^2}{(n + 1) H - 2 \kappa_i}  \\
& + \frac{H - \kappa_i}{3 H - 2 \kappa_j - \kappa_i} \frac{\beta}{\nu^{n + 1}} \frac{u_i}{u} (\nu^{n+1} - \kappa_i) H, \qquad  j \neq i.
\end{aligned}
\end{equation}

Now taking \eqref{eq3-14} and \eqref{eq3-17} into $(c)$ in \eqref{eq3-8} yields

\begin{equation} \label{eq3-18}
\begin{aligned}
&  \frac{\mu_1}{2} \frac{1}{(n + 1) H - 2 \kappa_i} \Bigg( (H - n \kappa_i) \\
& + \sum\limits_{j \neq i} \frac{H - \kappa_j}{3 H - 2 \kappa_j - \kappa_i} \frac{H \kappa_i - (n + 2) H \kappa_j - 2 \kappa_i^2 + (n + 2) \kappa_i \kappa_j + H^2}{H - \kappa_i} \Bigg) \\
= & \frac{\beta}{\nu^{n + 1}} \frac{u_i}{u} (\kappa_i - \nu^{n + 1}) H \sum\limits_{j \neq i} \frac{1}{3 H - 2 \kappa_j - \kappa_i} .
\end{aligned}
\end{equation}
We note that
\begin{equation} \label{eq3-19}
\begin{aligned}
& (H - n \kappa_i) + \sum\limits_{j \neq i} \frac{H - \kappa_j}{3 H - 2 \kappa_j - \kappa_i} \frac{H \kappa_i - (n + 2) H \kappa_j - 2 \kappa_i^2 + (n + 2) \kappa_i \kappa_j + H^2}{H - \kappa_i} \\
= &  (H - n \kappa_i) + \sum\limits_{j \neq i} \frac{H - \kappa_j}{3 H - 2 \kappa_j - \kappa_i} \Big( H + 2 \kappa_i - (n + 2) \kappa_j \Big) \\
= & - \frac{n^2 + n}{4} H + \frac{n^2 + 3 n - 2}{4} \kappa_i + \frac{1}{4} (H - \kappa_i) \Big( (3 n + 4) H - (n + 6) \kappa_i \Big) R_i,
\end{aligned}
\end{equation}
where
\[ R_i := \sum\limits_{j \neq i} \frac{1}{3 H - 2 \kappa_j - \kappa_i}. \]
Taking \eqref{eq3-19} into \eqref{eq3-18} we obtain
\begin{equation} \label{eq3-20}
\begin{aligned}
\frac{\mu_1}{2}
= \frac{ \frac{\beta}{\nu^{n + 1}} \frac{u_i}{u} (\kappa_i - \nu^{n + 1}) H \Big( (n + 1) H - 2 \kappa_i \Big) R_i }{ - \frac{n^2 + n}{4} H + \frac{n^2 + 3 n - 2}{4} \kappa_i + \frac{1}{4} (H - \kappa_i) \Big( (3 n + 4) H - (n + 6) \kappa_i \Big) R_i}.
\end{aligned}
\end{equation}
Taking \eqref{eq3-20} into \eqref{eq3-14} and \eqref{eq3-17} we obtain

\begin{equation*}
\begin{aligned}
t_i =
& \frac{\beta}{\nu^{n + 1}} \frac{u_i}{u} (\kappa_i - \nu^{n + 1}) H \frac{\frac{n^2 + n}{4} H - \frac{n^2 + 3 n - 2}{4} \kappa_i - \frac{1}{4} (H - \kappa_i) \Big( 3 n H + (3 n - 6) \kappa_i \Big) R_i }{ - \frac{n^2 + n}{4} H + \frac{n^2 + 3 n - 2}{4} \kappa_i + \frac{1}{4} (H - \kappa_i) \Big( (3 n + 4) H - (n + 6) \kappa_i \Big) R_i}.
\end{aligned}
\end{equation*}

\begin{footnotesize}
\begin{equation*}
\begin{aligned}
& t_j
= \frac{\beta}{\nu^{n + 1}} \frac{u_i}{u} (\kappa_i - \nu^{n + 1}) H \frac{H - \kappa_i}{3 H - 2 \kappa_j - \kappa_i} \\
& \cdot \frac{\frac{n^2 + n}{4} H - \frac{n^2 + 3 n - 2}{4} \kappa_i + \bigg( - \frac{3 n}{4} H^2 + \Big( n + \frac{9}{2} \Big) H \kappa_i - (n + 3) H \kappa_j - \Big( \frac{n}{4} + \frac{3}{2} \Big) \kappa_i^2 - 2 \kappa_i \kappa_j + (n + 2) \kappa_j^2 \bigg) R_i}{ - \frac{n^2 + n}{4} H + \frac{n^2 + 3 n - 2}{4} \kappa_i + \frac{1}{4} (H - \kappa_i) \Big( (3 n + 4) H - (n + 6) \kappa_i \Big) R_i} \\
= & \frac{\beta}{\nu^{n + 1}} \frac{u_i}{u} (\kappa_i - \nu^{n + 1}) H  \\
& \cdot \frac{\Big( \frac{n^2 + n}{4} H - \frac{n^2 + 3 n - 2}{4} \kappa_i \Big) \frac{H - \kappa_i}{3 H - 2 \kappa_j - \kappa_i} + (H - \kappa_i) \bigg( - \frac{n}{4} H - \frac{n + 2}{2} \kappa_j + \Big( \frac{n}{4} + \frac{3}{2} \Big) \kappa_i \bigg) R_i}{ - \frac{n^2 + n}{4} H + \frac{n^2 + 3 n - 2}{4} \kappa_i + \frac{1}{4} (H - \kappa_i) \Big( (3 n + 4) H - (n + 6) \kappa_i \Big) R_i},  \quad  j \neq i.
\end{aligned}
\end{equation*}
\end{footnotesize}

We summarize the following theorem.

\begin{thm} \label{OptimalConcavity}
The minimum of the quadratic form
\begin{equation*}
\begin{aligned}
Q_i
= &  \frac{2 \sigma}{n} \frac{1}{H - \kappa_i} \sum\limits_{j \neq i} \frac{1}{H - \kappa_j} t_j^2  - \frac{\sigma}{n} \sum\limits_{j \neq l} \sum\limits_{\substack{p < q \\ p, q \notin \{j, l\} }} \frac{1}{(H - \kappa_p)(H - \kappa_q)} t_j t_l \\
& - \frac{\sigma}{n} \sum\limits_j \bigg( 2 \sum\limits_{\substack{p < q \\ p, q \neq j }} \frac{1}{(H - \kappa_p)(H - \kappa_q)} -  \sum\limits_{p < q} \frac{1}{(H - \kappa_p)(H - \kappa_q)} \bigg) t_j^2 \\
& - \frac{\beta^2}{(\nu^{n + 1})^2} \frac{u_i^2}{u^2} (\nu^{n+1} - \kappa_i)^2 H^2 \frac{\sigma}{n} \sum\limits_{p < q} \frac{1}{(H - \kappa_p)(H - \kappa_q)}
\end{aligned}
\end{equation*}
is achieved at
\begin{equation*}
\begin{aligned}
t_i =
& \frac{\beta}{\nu^{n + 1}} \frac{u_i}{u} (\kappa_i - \nu^{n + 1}) H \frac{\frac{n^2 + n}{4} H - \frac{n^2 + 3 n - 2}{4} \kappa_i - \frac{1}{4} (H - \kappa_i) \Big( 3 n H + (3 n - 6) \kappa_i \Big) R_i }{ - \frac{n^2 + n}{4} H + \frac{n^2 + 3 n - 2}{4} \kappa_i + \frac{1}{4} (H - \kappa_i) \Big( (3 n + 4) H - (n + 6) \kappa_i \Big) R_i},
\end{aligned}
\end{equation*}
and for $j \neq i$,
\begin{equation*}
\begin{aligned}
t_j
= & \frac{\beta}{\nu^{n + 1}} \frac{u_i}{u} (\kappa_i - \nu^{n + 1}) H  \\
& \cdot \frac{\Big( \frac{n^2 + n}{4} H - \frac{n^2 + 3 n - 2}{4} \kappa_i \Big) \frac{H - \kappa_i}{3 H - 2 \kappa_j - \kappa_i} + (H - \kappa_i) \bigg( - \frac{n}{4} H - \frac{n + 2}{2} \kappa_j + \Big( \frac{n}{4} + \frac{3}{2} \Big) \kappa_i \bigg) R_i}{ - \frac{n^2 + n}{4} H + \frac{n^2 + 3 n - 2}{4} \kappa_i + \frac{1}{4} (H - \kappa_i) \Big( (3 n + 4) H - (n + 6) \kappa_i \Big) R_i},
\end{aligned}
\end{equation*}
where
\[ R_i := \sum\limits_{j \neq i} \frac{1}{3 H - 2 \kappa_j - \kappa_i}. \]
\end{thm}

\begin{rem} \label{n=3}
For $n = 3$, $i = 3$, by Theorem \ref{OptimalConcavity}, we find that
\begin{equation*}
\begin{aligned}
t_3 = \frac{\beta}{\nu^{n + 1}} \frac{u_3}{u} (\kappa_3 - \nu^{n + 1}) H \frac{3 \kappa_1^2 \kappa_2 + 3 \kappa_1 \kappa_2^2 - \kappa_1 \kappa_2 \kappa_3 - 2 \kappa_1 \kappa_3^2 -
   2 \kappa_2 \kappa_3^2 - \kappa_3^3}{\kappa_1^3 + \kappa_2^3 - \kappa_1^2 \kappa_3 - \kappa_1 \kappa_2 \kappa_3 - \kappa_2^2 \kappa_3 + \kappa_3^3},
\end{aligned}
\end{equation*}
\begin{equation*}
\begin{aligned}
t_1
= & - \frac{\beta}{\nu^{n + 1}} \frac{u_3}{u} (\kappa_3 - \nu^{n + 1}) H \frac{(\kappa_1 + \kappa_2) (\kappa_1^2 + \kappa_1 \kappa_2 + \kappa_1 \kappa_3 - 2 \kappa_2 \kappa_3 - \kappa_3^2)}{
 \kappa_1^3 + \kappa_2^3 - \kappa_1^2 \kappa_3 - \kappa_1 \kappa_2 \kappa_3 - \kappa_2^2 \kappa_3 + \kappa_3^3},
\end{aligned}
\end{equation*}
\begin{equation*}
\begin{aligned}
t_2
= & - \frac{\beta}{\nu^{n + 1}} \frac{u_3}{u} (\kappa_3 - \nu^{n + 1}) H  \frac{(\kappa_1 + \kappa_2) (\kappa_1 \kappa_2 + \kappa_2^2 - 2 \kappa_1 \kappa_3 + \kappa_2 \kappa_3 - \kappa_3^2)}{
 \kappa_1^3 + \kappa_2^3 - \kappa_1^2 \kappa_3 - \kappa_1 \kappa_2 \kappa_3 - \kappa_2^2 \kappa_3 + \kappa_3^3}.
\end{aligned}
\end{equation*}
\begin{equation*}
\begin{aligned}
& 2 \sum\limits_{j \neq 3} \frac{f_3 - f_j}{\kappa_j - \kappa_3} h_{3jj}^2 = \frac{2 \sigma}{3} \frac{\beta^2}{(\nu^{n + 1})^2} \frac{u_3^2}{u^2} (\kappa_3 - \nu^{n + 1})^2 H^2 \\
& \cdot \frac{(\kappa_1 + \kappa_2) \Big(  \substack{\kappa_1^5 + \kappa_1^4 (2 \kappa_2 + 3 \kappa_3) + \kappa_1^3 (\kappa_2^2 + \kappa_3^2) + \kappa_1^2 (\kappa_2^3 - 6 \kappa_2^2 \kappa_3 - 8 \kappa_2 \kappa_3^2 + \kappa_3^3) \\
+ \kappa_1 (2 \kappa_2^4 - 8 \kappa_2^2 \kappa_3^2 - 4 \kappa_2 \kappa_3^3 + 3 \kappa_3^4) + \kappa_2^5 + 3 \kappa_2^4 \kappa_3 + \kappa_2^3 \kappa_3^2 + \kappa_2^2 \kappa_3^3 +
3 \kappa_2 \kappa_3^4 + 2 \kappa_3^5}  \Big)}{(\kappa_1 + \kappa_3) (\kappa_2 + \kappa_3) (\kappa_1^3 + \kappa_2^3 - \kappa_1^2 \kappa_3 - \kappa_1 \kappa_2 \kappa_3 - \kappa_2^2 \kappa_3 + \kappa_3^3)^2},
\end{aligned}
\end{equation*}
\begin{equation*}
\begin{aligned}
& - \sum\limits_{j l} f_{jl}  h_{jj3} h_{ll3} = \frac{\sigma}{3} \frac{\beta^2}{(\nu^{n + 1})^2} \frac{u_3^2}{u^2} (\kappa_3 - \nu^{n + 1})^2 H^2 \\
& \cdot \frac{\Big(\substack{6 \kappa_1^4 + 6 \kappa_1^3 (\kappa_2 - \kappa_3) - 2 \kappa_1^2 \kappa_3 (3 \kappa_2 + 2 \kappa_3) + 2 \kappa_1 (3 \kappa_2^3 - 3 \kappa_2^2 \kappa_3 - 5 \kappa_2 \kappa_3^2 + 3 \kappa_3^3) \\ + 6 \kappa_2^4 - 6 \kappa_2^3 \kappa_3 - 4 \kappa_2^2 \kappa_3^2 + 6 \kappa_2 \kappa_3^3 + 6 \kappa_3^4 } \Big)}{(\kappa_1^3 + \kappa_2^3 - \kappa_1^2 \kappa_3 - \kappa_1 \kappa_2 \kappa_3 - \kappa_2^2 \kappa_3 + \kappa_3^3)^2}.
\end{aligned}
\end{equation*}
Therefore,
\begin{equation*}
\begin{aligned}
Q_3 = & \frac{2 \sigma}{3} \frac{\beta^2}{(\nu^{n + 1})^2} \frac{u_3^2}{u^2} (\kappa_3 - \nu^{n + 1})^2 H^2 \\
& \cdot \frac{\Big( \substack{\kappa_1^3 + 6 \kappa_1^2 \kappa_2 + 7 \kappa_1^2 \kappa_3 + 6 \kappa_1 \kappa_2^2  + 13 \kappa_1 \kappa_2 \kappa_3 + 8 \kappa_1 \kappa_3^2 \\ + \kappa_2^3 + 7 \kappa_2^2 \kappa_3  + 8 \kappa_2 \kappa_3^2 + 3 \kappa_3^3 } \Big)}{(\kappa_1 + \kappa_3) (\kappa_2 + \kappa_3) (\kappa_1^3 + \kappa_2^3 - \kappa_1^2 \kappa_3 - \kappa_1 \kappa_2 \kappa_3 - \kappa_2^2 \kappa_3 + \kappa_3^3)}.
\end{aligned}
\end{equation*}
Similarly, we can compute $Q_1$ and $Q_2$. By this way, we can recover formula (4.33) in \cite{CSS2025}.
\end{rem}

By Theorem \ref{OptimalConcavity}, \eqref{eq4-9} reduces to
\begin{equation} \label{eq5-2}
\begin{aligned}
& \bigg( A - \frac{\beta \big( 1 + (\nu^{n+1})^2 \big)}{\nu^{n + 1}} H + n \bigg) \sigma + ( \beta - 1 ) H \bigg( \sum f_i + \sum f_i \kappa_i^2 \bigg) \\
& + \frac{2 \sigma}{n} \frac{1}{H - \kappa_i} \sum\limits_{j \neq i} \frac{1}{H - \kappa_j} t_j^2  - \frac{\sigma}{n} \sum\limits_{j \neq l} \sum\limits_{\substack{p < q \\ p, q \notin \{j, l\} }} \frac{1}{(H - \kappa_p)(H - \kappa_q)} t_j t_l \\
& - \frac{\sigma}{n} \sum\limits_j \bigg( 2 \sum\limits_{\substack{p < q \\ p, q \neq j }} \frac{1}{(H - \kappa_p)(H - \kappa_q)} -  \sum\limits_{p < q} \frac{1}{(H - \kappa_p)(H - \kappa_q)} \bigg) t_j^2 \\
& - \frac{\beta^2}{(\nu^{n + 1})^2} \frac{u_i^2}{u^2} (\nu^{n+1} - \kappa_i)^2 H^2 \frac{\sigma}{n} \sum\limits_{p < q} \frac{1}{(H - \kappa_p)(H - \kappa_q)}
\\
& + \frac{2 \beta}{\nu^{n + 1}} H \sum f_i \frac{u_i^2}{u^2} \big( \kappa_i - \nu^{n+1} \big) - \frac{\beta (\beta - 1)}{(\nu^{n + 1})^2} H \sum f_i  \frac{u_i^2}{u^2} (\nu^{n + 1} - \kappa_i)^2 \leq 0.
\end{aligned}
\end{equation}

Hereafter, we focus on dimension $n = 4$.  We shall analyze each $i$ level, where $i = 1, \ldots, 4$.

\vspace{4mm}

\section{For $i = 1$}

\vspace{4mm}

In this section, we focus on $i = 1$ terms. By Mathematica, we have

\begin{equation*}
\begin{aligned}
t_1 =
& - \frac{\beta}{\nu^{4 + 1}} \frac{u_1}{u} (\kappa_1 - \nu^{4 + 1}) H \\
& \cdot \frac{\Bigg( \substack{12 \kappa_1^4 + 56 \kappa_1^3 (\kappa_2 + \kappa_3 + \kappa_4) +
     \kappa_1^2 \big((\kappa_2 + \kappa_3 + \kappa_4) (67 \kappa_2 + 67 \kappa_3 + 67 \kappa_4) + 12 (\kappa_2 \kappa_3 + \kappa_2 \kappa_4 + \kappa_3 \kappa_4)\big) \\
     + k1 \big( (\kappa_2 + \kappa_3 + \kappa_4) (15 \kappa_2^2 + 15 \kappa_3^2 + 15 \kappa_4^2 + 26 \kappa_2 \kappa_3 + 26 \kappa_2 \kappa_4 + 26 \kappa_3 \kappa_4) - 12 \kappa_2 \kappa_3 \kappa_4 \big) \\ - 8 (\kappa_2 + \kappa_3 + \kappa_4) (6 \kappa_2^2 \kappa_3 + 6 \kappa_2 \kappa_3^2 + 6 \kappa_2^2 \kappa_4 + 13 \kappa_2 \kappa_3 \kappa_4 + 6 \kappa_3^2 \kappa_4 + 6 \kappa_2 \kappa_4^2 + 6 \kappa_3 \kappa_4^2)} \Bigg)}{\Bigg( \substack{ 12 \kappa_1^4 + 20 \kappa_1^3 (\kappa_2 + \kappa_3 + \kappa_4) +
     \kappa_1^2 \big( 12 (\kappa_2 \kappa_3 + \kappa_2 \kappa_4 + \kappa_3 \kappa_4) - 5 (\kappa_2 + \kappa_3 + \kappa_4)^2 \big) \\ - 2 \kappa_1 \big((\kappa_2 + \kappa_3 + \kappa_4) (\kappa_2^2 + \kappa_3^2 + \kappa_4^2 + 10 \kappa_2 \kappa_3 + 10 \kappa_2 \kappa_4 + 10 \kappa_3 \kappa_4) + 6 \kappa_2 \kappa_3 \kappa_4 \big) \\
    + (\kappa_2 + \kappa_3 + \kappa_4) (15 \kappa_2^3 + \kappa_2^2 \kappa_3 + \kappa_2 \kappa_3^2 + 15 \kappa_3^3 + \kappa_2^2 \kappa_4 - 2 \kappa_2 \kappa_3 \kappa_4 + \kappa_3^2 \kappa_4 + \kappa_2 \kappa_4^2 + \kappa_3 \kappa_4^2 + 15 \kappa_4^3)} \Bigg)},
\end{aligned}
\end{equation*}

\begin{equation*}
\begin{aligned}
t_2
= & - \frac{\beta}{\nu^{4 + 1}} \frac{u_1}{u} (\kappa_1 - \nu^{4 + 1}) H  \\
& \cdot \frac{(\kappa_2 + \kappa_3 + \kappa_4) \Bigg( \substack{ -12 \kappa_1^3 + 4 \kappa_1^2 \kappa_2 + 43 \kappa_1 \kappa_2^2 + 15 \kappa_2^3 - 38 \kappa_1^2 \kappa_3 + 7 \kappa_1 \kappa_2 \kappa_3 \\
+ 46 \kappa_2^2 \kappa_3 - 30 \kappa_1 \kappa_3^2 + 19 \kappa_2 \kappa_3^2 - 38 \kappa_1^2 \kappa_4 + 7 \kappa_1 \kappa_2 \kappa_4 + 46 \kappa_2^2 \kappa_4 \\
- 60 \kappa_1 \kappa_3 \kappa_4 + 34 \kappa_2 \kappa_3 \kappa_4 - 16 \kappa_3^2 \kappa_4 - 30 \kappa_1 \kappa_4^2 + 19 \kappa_2 \kappa_4^2 - 16 \kappa_3 \kappa_4^2 } \Bigg) }{\Bigg( \substack{ 12 \kappa_1^4 + 20 \kappa_1^3 (\kappa_2 + \kappa_3 + \kappa_4) +
     \kappa_1^2 \big( 12 (\kappa_2 \kappa_3 + \kappa_2 \kappa_4 + \kappa_3 \kappa_4) - 5 (\kappa_2 + \kappa_3 + \kappa_4)^2 \big) \\ - 2 \kappa_1 \big((\kappa_2 + \kappa_3 + \kappa_4) (\kappa_2^2 + \kappa_3^2 + \kappa_4^2 + 10 \kappa_2 \kappa_3 + 10 \kappa_2 \kappa_4 + 10 \kappa_3 \kappa_4) + 6 \kappa_2 \kappa_3 \kappa_4 \big) \\
    + (\kappa_2 + \kappa_3 + \kappa_4) (15 \kappa_2^3 + \kappa_2^2 \kappa_3 + \kappa_2 \kappa_3^2 + 15 \kappa_3^3 + \kappa_2^2 \kappa_4 - 2 \kappa_2 \kappa_3 \kappa_4 + \kappa_3^2 \kappa_4 + \kappa_2 \kappa_4^2 + \kappa_3 \kappa_4^2 + 15 \kappa_4^3)} \Bigg)},
\end{aligned}
\end{equation*}

\begin{equation*}
\begin{aligned}
t_3
= & - \frac{\beta}{\nu^{4 + 1}} \frac{u_1}{u} (\kappa_1 - \nu^{4 + 1}) H  \\
& \cdot \frac{(\kappa_2 + \kappa_3 + \kappa_4) \Bigg( \substack{ -12 \kappa_1^3 - 38 \kappa_1^2 \kappa_2 - 30 \kappa_1 \kappa_2^2 + 4 \kappa_1^2 \kappa_3
+ 7 \kappa_1 \kappa_2 \kappa_3 + 19 \kappa_2^2 \kappa_3 \\
+ 43 \kappa_1 \kappa_3^2 + 46 \kappa_2 \kappa_3^2 + 15 \kappa_3^3 - 38 \kappa_1^2 \kappa_4 - 60 \kappa_1 \kappa_2 \kappa_4 - 16 \kappa_2^2 \kappa_4 \\ + 7 \kappa_1 \kappa_3 \kappa_4
+ 34 \kappa_2 \kappa_3 \kappa_4 + 46 \kappa_3^2 \kappa_4 - 30 \kappa_1 \kappa_4^2 - 16 \kappa_2 \kappa_4^2 + 19 \kappa_3 \kappa_4^2 } \Bigg)}{\Bigg( \substack{ 12 \kappa_1^4 + 20 \kappa_1^3 (\kappa_2 + \kappa_3 + \kappa_4) +
     \kappa_1^2 \big( 12 (\kappa_2 \kappa_3 + \kappa_2 \kappa_4 + \kappa_3 \kappa_4) - 5 (\kappa_2 + \kappa_3 + \kappa_4)^2 \big) \\ - 2 \kappa_1 \big((\kappa_2 + \kappa_3 + \kappa_4) (\kappa_2^2 + \kappa_3^2 + \kappa_4^2 + 10 \kappa_2 \kappa_3 + 10 \kappa_2 \kappa_4 + 10 \kappa_3 \kappa_4) + 6 \kappa_2 \kappa_3 \kappa_4 \big) \\
    + (\kappa_2 + \kappa_3 + \kappa_4) (15 \kappa_2^3 + \kappa_2^2 \kappa_3 + \kappa_2 \kappa_3^2 + 15 \kappa_3^3 + \kappa_2^2 \kappa_4 - 2 \kappa_2 \kappa_3 \kappa_4 + \kappa_3^2 \kappa_4 + \kappa_2 \kappa_4^2 + \kappa_3 \kappa_4^2 + 15 \kappa_4^3)} \Bigg)},
\end{aligned}
\end{equation*}

\begin{equation*}
\begin{aligned}
t_4
= & \frac{\beta}{\nu^{4 + 1}} \frac{u_1}{u} (\kappa_1 - \nu^{4 + 1}) H  \\
& \cdot \frac{(\kappa_2 + \kappa_3 + \kappa_4) \Bigg( \substack{12 \kappa_1^3 + 38 \kappa_1^2 \kappa_2 + 30 \kappa_1 \kappa_2^2 + 38 \kappa_1^2 \kappa_3 + 60 \kappa_1 \kappa_2 \kappa_3 + 16 \kappa_2^2 \kappa_3 \\
+ 30 \kappa_1 \kappa_3^2 + 16 \kappa_2 \kappa_3^2 - 4 \kappa_1^2 \kappa_4 - 7 \kappa_1 \kappa_2 \kappa_4 - 19 \kappa_2^2 \kappa_4 - 7 \kappa_1 \kappa_3 \kappa_4 \\ - 34 \kappa_2 \kappa_3 \kappa_4 - 19 \kappa_3^2 \kappa_4 - 43 \kappa_1 \kappa_4^2 - 46 \kappa_2 \kappa_4^2 - 46 \kappa_3 \kappa_4^2 - 15 \kappa_4^3} \Bigg)}{\Bigg( \substack{ 12 \kappa_1^4 + 20 \kappa_1^3 (\kappa_2 + \kappa_3 + \kappa_4) +
     \kappa_1^2 \big( 12 (\kappa_2 \kappa_3 + \kappa_2 \kappa_4 + \kappa_3 \kappa_4) - 5 (\kappa_2 + \kappa_3 + \kappa_4)^2 \big) \\ - 2 \kappa_1 \big((\kappa_2 + \kappa_3 + \kappa_4) (\kappa_2^2 + \kappa_3^2 + \kappa_4^2 + 10 \kappa_2 \kappa_3 + 10 \kappa_2 \kappa_4 + 10 \kappa_3 \kappa_4) + 6 \kappa_2 \kappa_3 \kappa_4 \big) \\
    + (\kappa_2 + \kappa_3 + \kappa_4) (15 \kappa_2^3 + \kappa_2^2 \kappa_3 + \kappa_2 \kappa_3^2 + 15 \kappa_3^3 + \kappa_2^2 \kappa_4 - 2 \kappa_2 \kappa_3 \kappa_4 + \kappa_3^2 \kappa_4 + \kappa_2 \kappa_4^2 + \kappa_3 \kappa_4^2 + 15 \kappa_4^3)} \Bigg)}.
\end{aligned}
\end{equation*}

We note that
\begin{footnotesize}
\begin{equation*}
\begin{aligned}
& - \sum\limits_{j l} f_{jl}  h_{jj1} h_{ll1} \\
= & - \frac{\sigma}{4} \frac{1}{(H - \kappa_3)(H - \kappa_4)} (t_1 + t_2)^2 - \frac{\sigma}{4} \frac{1}{(H - \kappa_2)(H - \kappa_4)} (t_1 + t_3)^2 \\
  & - \frac{\sigma}{4} \frac{1}{(H - \kappa_2)(H - \kappa_3)} (t_1 + t_4)^2  - \frac{\sigma}{4} \frac{1}{(H - \kappa_1)(H - \kappa_4)} (t_2 + t_3)^2 \\
    & - \frac{\sigma}{4} \frac{1}{(H - \kappa_1)(H - \kappa_3)} (t_2 + t_4)^2  - \frac{\sigma}{4} \frac{1}{(H - \kappa_1)(H - \kappa_2)} (t_3 + t_4)^2 \\
& + \frac{\sigma}{4} \bigg( \frac{1}{(H - \kappa_1)(H - \kappa_2)} + \frac{1}{(H - \kappa_1)(H - \kappa_3)} + \frac{1}{(H - \kappa_1)(H - \kappa_4)} \bigg) t_1^2 \\
& + \frac{\sigma}{4} \bigg( \frac{1}{(H - \kappa_1)(H - \kappa_2)} + \frac{1}{(H - \kappa_2)(H - \kappa_3)}  + \frac{1}{(H - \kappa_2)(H - \kappa_4)}  \bigg) t_2^2 \\
& + \frac{\sigma}{4} \bigg( \frac{1}{(H - \kappa_1)(H - \kappa_3)} + \frac{1}{(H - \kappa_2)(H - \kappa_3)}  + \frac{1}{(H - \kappa_3)(H - \kappa_4)}  \bigg) t_3^2 \\
& + \frac{\sigma}{4} \bigg( \frac{1}{(H - \kappa_1)(H - \kappa_4)} + \frac{1}{(H - \kappa_2)(H - \kappa_4)}  + \frac{1}{(H - \kappa_3)(H - \kappa_4)} \bigg) t_4^2 \\
& - \frac{\beta^2}{(\nu^{4 + 1})^2} \frac{u_1^2}{u^2} (\nu^{4 + 1} - \kappa_1)^2 H^2 \frac{\sigma}{4} \sum\limits_{p < q} \frac{1}{(H - \kappa_p)(H - \kappa_q)} \\
= & \frac{\sigma}{4} \Bigg(  \frac{1}{(H - \kappa_3)(H - \kappa_4)} \bigg( t_3^2 + t_4^2 - (t_1 + t_2)^2 - \frac{\beta^2}{(\nu^{4 + 1})^2} \frac{u_1^2} {u^2} (\nu^{4 + 1} - \kappa_1)^2 H^2 \bigg) \\
& + \frac{1}{(H - \kappa_2)(H - \kappa_4)} \bigg( t_2^2 + t_4^2 - (t_1 + t_3)^2 - \frac{\beta^2}{(\nu^{4 + 1})^2} \frac{u_1^2} {u^2} (\nu^{4 + 1} - \kappa_1)^2 H^2 \bigg) \\
& + \frac{1}{(H - \kappa_2)(H - \kappa_3)} \bigg( t_2^2 + t_3^2 - (t_1 + t_4)^2 - \frac{\beta^2}{(\nu^{4 + 1})^2} \frac{u_1^2}{u^2} (\nu^{4 + 1} - \kappa_1)^2 H^2 \bigg) \\
& + \frac{1}{(H - \kappa_1) (H - \kappa_4)} \bigg( t_1^2 + t_4^2 - (t_2 + t_3)^2 - \frac{\beta^2}{(\nu^{4 + 1})^2} \frac{u_1^2}{u^2} (\nu^{4 + 1} - \kappa_1)^2 H^2 \bigg) \\
& + \frac{1}{(H - \kappa_1)(H - \kappa_3)} \bigg( t_1^2 + t_3^2 - (t_2 + t_4)^2 - \frac{\beta^2}{(\nu^{4 + 1})^2} \frac{u_1^2}{u^2} (\nu^{4 + 1} - \kappa_1)^2 H^2 \bigg) \\
& + \frac{1}{(H - \kappa_1)(H - \kappa_2)} \bigg( t_1^2 + t_2^2 - (t_3 + t_4)^2 - \frac{\beta^2}{(\nu^{4 + 1})^2} \frac{u_1^2}{u^2} (\nu^{4 + 1} - \kappa_1)^2 H^2 \bigg) \Bigg).
\end{aligned}
\end{equation*}
\end{footnotesize}

By Mathematica, we compute that
\begin{footnotesize}
\begin{equation} \label{eq5-3}
\begin{aligned}
& - \sum\limits_{j l} f_{jl}  h_{jj1} h_{ll1} \\
= &  \frac{4 \sigma}{4} \frac{\beta^2}{(\nu^{4 + 1})^2} \frac{u_1^2}{u^2} (\kappa_1 - \nu^{4 + 1})^2 H^2  \\
& \cdot \frac{\left( \substack{ 432 \kappa_1^6 + 1728 \kappa_1^5 (\kappa_2 + \kappa_3 + \kappa_4) +
     216 \kappa_1^4 \big( 10 (\kappa_2 + \kappa_3 + \kappa_4)^2 + \kappa_2 \kappa_3 + \kappa_2 \kappa_4 + \kappa_3 \kappa_4 \big) \\
+ 8 k1^3 \big( (\kappa_2 + \kappa_3 + \kappa_4) (46 \kappa_2^2 + 46 \kappa_3^2 + 46 \kappa_4^2 + 182 \kappa_2 \kappa_3 + 182 \kappa_2 \kappa_4 + 182 \kappa_3 \kappa_4) + 27 \kappa_2 \kappa_3 \kappa_4 \big) \\
 - 3 \kappa_1^2 \big( (\kappa_2 + \kappa_3 + \kappa_4) (261 \kappa_2^3 + 261 \kappa_3^3 + 261 \kappa_4^3 + 601 \kappa_2^2 \kappa_3 + 601 \kappa_2 \kappa_3^2 + 601 \kappa_2^2 \kappa_4 \\
 + 601 \kappa_2 \kappa_4^2 + 601 \kappa_3^2 \kappa_4 + 601 \kappa_3 \kappa_4^2 + 1368 \kappa_2 \kappa_3 \kappa_4) - 24 (\kappa_2^2 \kappa_3^2 + \kappa_2^2 \kappa_4^2 + \kappa_3^2 \kappa_4^2) \big) \\
 + 6 \kappa_1 (\kappa_2 + \kappa_3 + \kappa_4) (45 \kappa_2^4 - 72 \kappa_2^3 \kappa_3 - 266 \kappa_2^2 \kappa_3^2 - 72 \kappa_2 \kappa_3^3 + 45 \kappa_3^4 - 72 \kappa_2^3 \kappa_4
 - 664 \kappa_2^2 \kappa_3 \kappa_4 \\
 - 664 \kappa_2 \kappa_3^2 \kappa_4 - 72 \kappa_3^3 \kappa_4 - 266 \kappa_2^2 \kappa_4^2 - 664 \kappa_2 \kappa_3 \kappa_4^2 - 266 \kappa_3^2 \kappa_4^2 - 72 \kappa_2 \kappa_4^3 - 72 \kappa_3 \kappa_4^3 + 45 \kappa_4^4) \\
+ (\kappa_2 + \kappa_3 + \kappa_4)^2 \big(675 \kappa_2^4 + 675 \kappa_3^4 + 540 \kappa_3^3 \kappa_4 + 2 \kappa_3^2 \kappa_4^2 + 540 \kappa_3 \kappa_4^3 + 675 \kappa_4^4 \\ + 540 \kappa_2^3 (\kappa_3 + \kappa_4) +
 2 \kappa_2^2 (\kappa_3^2 + 322 \kappa_3 \kappa_4 + \kappa_4^2) + \kappa_2 (540 \kappa_3^3 + 644 \kappa_3^2 \kappa_4 + 644 \kappa_3 \kappa_4^2 + 540 \kappa_4^3) \big) } \right)}{\left( \substack{ 12 \kappa_1^4 + 20 \kappa_1^3 (\kappa_2 + \kappa_3 + \kappa_4) +
 \kappa_1^2 \big( 12 (\kappa_2 \kappa_3 + \kappa_2 \kappa_4 + \kappa_3 \kappa_4) - 5 (\kappa_2 + \kappa_3 + \kappa_4)^2 \big) \\
 - 2 \kappa_1 \big( (\kappa_2 + \kappa_3 + \kappa_4) (\kappa_2^2 + \kappa_3^2 + \kappa_4^2 + 10 \kappa_2 \kappa_3 + 10 \kappa_2 \kappa_4 + 10 \kappa_3 \kappa_4) + 6 \kappa_2 \kappa_3 \kappa_4 \big) \\
  + (\kappa_2 + \kappa_3 + \kappa_4) (15 \kappa_2^3 + \kappa_2^2 \kappa_3 + \kappa_2 \kappa_3^2 + 15 \kappa_3^3 + \kappa_2^2 \kappa_4 - 2 \kappa_2 \kappa_3 \kappa_4 + \kappa_3^2 \kappa_4 + \kappa_2 \kappa_4^2 + \kappa_3 \kappa_4^2 + 15 \kappa_4^3) } \right)^2}.
\end{aligned}
\end{equation}
\end{footnotesize}

\begin{footnotesize}
\begin{equation*}
\begin{aligned}
& \frac{2 \sigma}{4} \frac{1}{H - \kappa_1}  \frac{1}{H - \kappa_2} t_2^2 \\
= & \frac{2 \sigma}{4}  \frac{\beta^2}{(\nu^{4 + 1})^2} \frac{u_1^2}{u^2} (\kappa_1 - \nu^{4 + 1})^2 H^2 \\
& \cdot \frac{(\kappa_2 + \kappa_3 + \kappa_4) \Bigg( \substack{ 12 \kappa_1^3 - 4 \kappa_1^2 \kappa_2 - 43 \kappa_1 \kappa_2^2 - 15 \kappa_2^3
+ 38 \kappa_1^2 \kappa_3 - 7 \kappa_1 \kappa_2 \kappa_3 \\
- 46 \kappa_2^2 \kappa_3 + 30 \kappa_1 \kappa_3^2 - 19 \kappa_2 \kappa_3^2 + 38 \kappa_1^2 \kappa_4
- 7 \kappa_1 \kappa_2 \kappa_4 - 46 \kappa_2^2 \kappa_4 + 60 \kappa_1 \kappa_3 \kappa_4 \\
- 34 \kappa_2 \kappa_3 \kappa_4 + 16 \kappa_3^2 \kappa_4 + 30 \kappa_1 \kappa_4^2 - 19 \kappa_2 \kappa_4^2 + 16 \kappa_3 \kappa_4^2 } \Bigg)^2}{(\kappa_1 + \kappa_3 + \kappa_4) \left( \substack{ 12 \kappa_1^4 + 20 \kappa_1^3 (\kappa_2 + \kappa_3 + \kappa_4) +
 \kappa_1^2 \big( 12 (\kappa_2 \kappa_3 + \kappa_2 \kappa_4 + \kappa_3 \kappa_4) - 5 (\kappa_2 + \kappa_3 + \kappa_4)^2 \big) \\
 - 2 \kappa_1 \big( (\kappa_2 + \kappa_3 + \kappa_4) (\kappa_2^2 + \kappa_3^2 + \kappa_4^2 + 10 \kappa_2 \kappa_3 + 10 \kappa_2 \kappa_4 + 10 \kappa_3 \kappa_4) + 6 \kappa_2 \kappa_3 \kappa_4 \big) \\
  + (\kappa_2 + \kappa_3 + \kappa_4) (15 \kappa_2^3 + \kappa_2^2 \kappa_3 + \kappa_2 \kappa_3^2 + 15 \kappa_3^3 + \kappa_2^2 \kappa_4 - 2 \kappa_2 \kappa_3 \kappa_4 + \kappa_3^2 \kappa_4 + \kappa_2 \kappa_4^2 + \kappa_3 \kappa_4^2 + 15 \kappa_4^3) } \right)^2}.
\end{aligned}
\end{equation*}
\end{footnotesize}

\begin{footnotesize}
\begin{equation*}
\begin{aligned}
& \frac{2 \sigma}{4} \frac{1}{H - \kappa_1} \frac{1}{H - \kappa_3} t_3^2 \\
= &  \frac{2 \sigma}{4}  \frac{\beta^2}{(\nu^{4 + 1})^2} \frac{u_1^2}{u^2} (\kappa_1 - \nu^{4 + 1})^2 H^2 \\
& \cdot \frac{(\kappa_2 + \kappa_3 + \kappa_4) \Bigg( \substack{ 12 \kappa_1^3 + 38 \kappa_1^2 \kappa_2 + 30 \kappa_1 \kappa_2^2 - 4 \kappa_1^2 \kappa_3 - 7 \kappa_1 \kappa_2 \kappa_3 \\
- 19 \kappa_2^2 \kappa_3 - 43 \kappa_1 \kappa_3^2 - 46 \kappa_2 \kappa_3^2 - 15 \kappa_3^3 +
38 \kappa_1^2 \kappa_4 + 60 \kappa_1 \kappa_2 \kappa_4 \\
+ 16 \kappa_2^2 \kappa_4 - 7 \kappa_1 \kappa_3 \kappa_4 - 34 \kappa_2 \kappa_3 \kappa_4 - 46 \kappa_3^2 \kappa_4 + 30 \kappa_1 \kappa_4^2 + 16 \kappa_2 \kappa_4^2 - 19 \kappa_3 \kappa_4^2 } \Bigg)^2}{(\kappa_1 + \kappa_2 + \kappa_4) \left( \substack{ 12 \kappa_1^4 + 20 \kappa_1^3 (\kappa_2 + \kappa_3 + \kappa_4) +
 \kappa_1^2 \big( 12 (\kappa_2 \kappa_3 + \kappa_2 \kappa_4 + \kappa_3 \kappa_4) - 5 (\kappa_2 + \kappa_3 + \kappa_4)^2 \big) \\
 - 2 \kappa_1 \big( (\kappa_2 + \kappa_3 + \kappa_4) (\kappa_2^2 + \kappa_3^2 + \kappa_4^2 + 10 \kappa_2 \kappa_3 + 10 \kappa_2 \kappa_4 + 10 \kappa_3 \kappa_4) + 6 \kappa_2 \kappa_3 \kappa_4 \big) \\
  + (\kappa_2 + \kappa_3 + \kappa_4) (15 \kappa_2^3 + \kappa_2^2 \kappa_3 + \kappa_2 \kappa_3^2 + 15 \kappa_3^3 + \kappa_2^2 \kappa_4 - 2 \kappa_2 \kappa_3 \kappa_4 + \kappa_3^2 \kappa_4 + \kappa_2 \kappa_4^2 + \kappa_3 \kappa_4^2 + 15 \kappa_4^3) } \right)^2}.
\end{aligned}
\end{equation*}
\end{footnotesize}

\begin{footnotesize}
\begin{equation*}
\begin{aligned}
& \frac{2 \sigma}{4} \frac{1}{H - \kappa_1} \frac{1}{H - \kappa_4} t_4^2 \\
= &  \frac{2 \sigma}{4}  \frac{\beta^2}{(\nu^{4 + 1})^2} \frac{u_1^2}{u^2} (\kappa_1 - \nu^{4 + 1})^2 H^2 \\
& \cdot \frac{(\kappa_2 + \kappa_3 + \kappa_4) \Bigg( \substack{ 12 \kappa_1^3 + 38 \kappa_1^2 \kappa_2 + 30 \kappa_1 \kappa_2^2 + 38 \kappa_1^2 \kappa_3 + 60 \kappa_1 \kappa_2 \kappa_3 \\
+ 16 \kappa_2^2 \kappa_3 + 30 \kappa_1 \kappa_3^2 + 16 \kappa_2 \kappa_3^2 - 4 \kappa_1^2 \kappa_4 -
 7 \kappa_1 \kappa_2 \kappa_4 - 19 \kappa_2^2 \kappa_4 - 7 \kappa_1 \kappa_3 \kappa_4 \\
- 34 \kappa_2 \kappa_3 \kappa_4 - 19 \kappa_3^2 \kappa_4 - 43 \kappa_1 \kappa_4^2 - 46 \kappa_2 \kappa_4^2 - 46 \kappa_3 \kappa_4^2 - 15 \kappa_4^3} \Bigg)^2}{(\kappa_1 + \kappa_2 + \kappa_3) \left( \substack{ 12 \kappa_1^4 + 20 \kappa_1^3 (\kappa_2 + \kappa_3 + \kappa_4) +
 \kappa_1^2 \big( 12 (\kappa_2 \kappa_3 + \kappa_2 \kappa_4 + \kappa_3 \kappa_4) - 5 (\kappa_2 + \kappa_3 + \kappa_4)^2 \big) \\
 - 2 \kappa_1 \big( (\kappa_2 + \kappa_3 + \kappa_4) (\kappa_2^2 + \kappa_3^2 + \kappa_4^2 + 10 \kappa_2 \kappa_3 + 10 \kappa_2 \kappa_4 + 10 \kappa_3 \kappa_4) + 6 \kappa_2 \kappa_3 \kappa_4 \big) \\
  + (\kappa_2 + \kappa_3 + \kappa_4) (15 \kappa_2^3 + \kappa_2^2 \kappa_3 + \kappa_2 \kappa_3^2 + 15 \kappa_3^3 + \kappa_2^2 \kappa_4 - 2 \kappa_2 \kappa_3 \kappa_4 + \kappa_3^2 \kappa_4 + \kappa_2 \kappa_4^2 + \kappa_3 \kappa_4^2 + 15 \kappa_4^3) } \right)^2}.
\end{aligned}
\end{equation*}
\end{footnotesize}

For $\kappa_1 \geq 1$, and any $\beta > 1$,

\begin{footnotesize}
\begin{equation*}
\begin{aligned}
& - \sum\limits_{j l} f_{jl}  h_{jj1} h_{ll1} +  \frac{2 \sigma}{4} \frac{1}{H - \kappa_1}  \frac{1}{H - \kappa_2} t_2^2 + \frac{2 \sigma}{4} \frac{1}{H - \kappa_1} \frac{1}{H - \kappa_3} t_3^2 \\
& + \frac{2 \sigma}{4} \frac{1}{H - \kappa_1} \frac{1}{H - \kappa_4} t_4^2
- \frac{\beta^2}{(\nu^{n + 1})^2} H  f_1  \frac{u_1^2}{u^2} (\kappa_1 - \nu^{4 + 1})^2 \\
= & \frac{\sigma}{4} \frac{\beta^2}{(\nu^{n + 1})^2} \frac{u_1^2}{u^2} (\nu^{4 + 1} - \kappa_1)^2 H \frac{1}{(\kappa_1 + \kappa_2 + \kappa_3) (\kappa_1 + \kappa_2 + \kappa_4) (\kappa_1 + \kappa_3 + \kappa_4)} \\
& \cdot \frac{\left( \substack{ 108 \kappa_1^6 + 732 \kappa_1^5 (\kappa_2 + \kappa_3 + \kappa_4)
+ k_1^4 \big( 1919 (\kappa_2 + \kappa_3 + \kappa_4)^2 + 24 (\kappa_2 \kappa_3 + \kappa_2 \kappa_4 + \kappa_3 \kappa_4) \big) \\
+ 4 \kappa_1^3 \big( 610 (\kappa_2 + \kappa_3 + \kappa_4)^3 + 88 (\kappa_2 + \kappa_3 + \kappa_4) (\kappa_2 \kappa_3 + \kappa_2 \kappa_4 + \kappa_3 \kappa_4) + 27 \kappa_2 \kappa_3 \kappa_4 \big) \\
+ 3 \kappa_1^2 \big( 506 (\kappa_2 + \kappa_3 + \kappa_4)^4 +
 317 (\kappa_2 + \kappa_3 + \kappa_4)^2 (\kappa_2 \kappa_3 + \kappa_2 \kappa_4 + \kappa_3 \kappa_4) \\
 + 4 (\kappa_2 \kappa_3 + \kappa_2 \kappa_4 + \kappa_3 \kappa_4)^2 - 58 (\kappa_2 + \kappa_3 + \kappa_4) \kappa_2 \kappa_3 \kappa_4 \big) \\
 + 2 \kappa_1 \big( 194 (\kappa_2 + \kappa_3 + \kappa_4)^5 + 449 (\kappa_2 + \kappa_3 + \kappa_4)^3 (\kappa_2 \kappa_3 + \kappa_2 \kappa_4 + \kappa_3 \kappa_4) -
 271 (\kappa_2 + \kappa_3 + \kappa_4)^2 \kappa_2 \kappa_3 \kappa_4 \\
 + 16 (\kappa_2 + \kappa_3 + \kappa_4) (\kappa_2 \kappa_3 + \kappa_2 \kappa_4 + \kappa_3 \kappa_4)^2 + 6 (\kappa_2 \kappa_3 + \kappa_2 \kappa_4 + \kappa_3 \kappa_4) \kappa_2 \kappa_3 \kappa_4 \big) \\
+ (\kappa_2 + \kappa_3 + \kappa_4) \big( 15 (\kappa_2 + \kappa_3 + \kappa_4)^5 + 271 (\kappa_2 + \kappa_3 + \kappa_4)^3 (\kappa_2 \kappa_3 + \kappa_2 \kappa_4 + \kappa_3 \kappa_4)
\\ + 36 (\kappa_2 + \kappa_3 + \kappa_4) (\kappa_2 \kappa_3 + \kappa_2 \kappa_4 + \kappa_3 \kappa_4)^2 -
 260 (\kappa_2 + \kappa_3 + \kappa_4)^2 \kappa_2 \kappa_3 \kappa_4 - 40 \kappa_2 \kappa_3 \kappa_4 (\kappa_2 \kappa_3 + \kappa_2 \kappa_4 + \kappa_3 \kappa_4) \big) } \right) }{\left( \substack{ 12 \kappa_1^4 + 20 \kappa_1^3 (\kappa_2 + \kappa_3 + \kappa_4) +
 \kappa_1^2 \big( 12 (\kappa_2 \kappa_3 + \kappa_2 \kappa_4 + \kappa_3 \kappa_4) - 5 (\kappa_2 + \kappa_3 + \kappa_4)^2 \big) \\
 - 2 \kappa_1 \big( (\kappa_2 + \kappa_3 + \kappa_4) (\kappa_2^2 + \kappa_3^2 + \kappa_4^2 + 10 \kappa_2 \kappa_3 + 10 \kappa_2 \kappa_4 + 10 \kappa_3 \kappa_4) + 6 \kappa_2 \kappa_3 \kappa_4 \big) \\
  + (\kappa_2 + \kappa_3 + \kappa_4) (15 \kappa_2^3 + \kappa_2^2 \kappa_3 + \kappa_2 \kappa_3^2 + 15 \kappa_3^3 + \kappa_2^2 \kappa_4 - 2 \kappa_2 \kappa_3 \kappa_4 + \kappa_3^2 \kappa_4 + \kappa_2 \kappa_4^2 + \kappa_3 \kappa_4^2 + 15 \kappa_4^3) } \right)}.
\end{aligned}
\end{equation*}
\end{footnotesize}

For convenience, we denote
\[ \sigma_1 = \kappa_2 + \kappa_3 + \kappa_4, \quad \sigma_2 = \kappa_2 \kappa_3 + \kappa_2 \kappa_4 + \kappa_3 \kappa_4, \quad \sigma_3 = \kappa_2 \kappa_3 \kappa_4. \]
Now we prove the nonnegativity of the top:
\begin{equation} \label{eq5-4}
\begin{aligned}
T_1 := & 108 \kappa_1^6 + 732 \kappa_1^5 \sigma_1
+ k_1^4 ( 1919 \sigma_1^2 + 24 \sigma_2 )
+ 4 \kappa_1^3 ( 610 \sigma_1^3 + 88 \sigma_1 \sigma_2 + 27 \sigma_3 ) \\
& + 3 \kappa_1^2 ( 506 \sigma_1^4 +
 317 \sigma_1^2 \sigma_2
 + 4 \sigma_2^2 - 58 \sigma_1 \sigma_3 ) \\
& + 2 \kappa_1 ( 194 \sigma_1^5 + 449 \sigma_1^3 \sigma_2 -
 271 \sigma_1^2 \sigma_3
 + 16 \sigma_1 \sigma_2^2 + 6 \sigma_2 \sigma_3 ) \\
& + \sigma_1 ( 15 \sigma_1^5 + 271 \sigma_1^3 \sigma_2
 + 36 \sigma_1 \sigma_2^2 -
 260 \sigma_1^2 \sigma_3 - 40 \sigma_3 \sigma_2 ) \\
= & 108 \kappa_1^6  + 24  k_1^4 \sigma_2  + 108 \kappa_1^3 \sigma_3 + 12 \kappa_1^2 \sigma_2^2 + 12 \kappa_1 \sigma_2 \sigma_3  \\
& + (732 \kappa_1^5 + 88 \cdot 4 \kappa_1^3 \sigma_2 - 58 \cdot 3 \kappa_1^2 \sigma_3  + 32 \kappa_1  \sigma_2^2  - 40 \sigma_3 \sigma_2 ) \sigma_1 \\
& + (1919 k_1^4 +
 317 \cdot 3 \kappa_1^2 \sigma_2  -
 271 \cdot 2 \kappa_1 \sigma_3 + 36 \sigma_2^2) \sigma_1^2  \\
& + (610 \cdot 4 \kappa_1^3 + 449 \cdot 2 \kappa_1 \sigma_2 -
 260 \sigma_3) \sigma_1^3  \\
& + ( 506 \cdot 3 \kappa_1^2 + 271 \sigma_2) \sigma_1^4 \\
& + 194 \cdot 2 \kappa_1  \sigma_1^5 \\
& + 15 \sigma_1^6 .
\end{aligned}
\end{equation}
We need the following propositions.
\begin{prop} \label{Prop0}
For admissible hypersurfaces, we have
\[ \kappa_1 \geq \kappa_2 > 0. \]
\[ \kappa_2 + \kappa_3 > 0. \]
Moreover, for admissible solutions, we have
\[ 0 < \kappa_2 + \kappa_3 + \kappa_4 \leq \frac{(3 \sigma)^{\frac{4}{3}}}{\kappa_1^{\frac{1}{3}}}. \]
\end{prop}
\begin{proof}
By admissibility of solutions, we know that
\[ \kappa_2 + \kappa_3 + \kappa_4 > 0, \]
and $\kappa_1 \geq \kappa_2 > 0$.
It follows that
\[ \kappa_2 + \kappa_3 > 0. \]
For otherwise, $\kappa_3 < 0$, but then $\kappa_4 \leq \kappa_3 < 0$ and consequently
\[ ( \kappa_2 + \kappa_3 ) + \kappa_4 < 0, \]
contradicting with admissibility.

Next, by equation
\[ (\kappa_1 + \kappa_2 + \kappa_3)(\kappa_1 + \kappa_2 + \kappa_4)(\kappa_1 + \kappa_3 + \kappa_4)(\kappa_2 + \kappa_3 + \kappa_4) = (3 \sigma)^4, \]
we can deduce that
\[ \kappa_1 (\kappa_2 + \kappa_3 + \kappa_4)^3 \leq (3 \sigma)^4. \]
\end{proof}

\begin{prop} \label{Prop1}
\[  \kappa_2 \kappa_3 + \kappa_2 \kappa_4 + \kappa_3 \kappa_4 + 3 \kappa_1^2 \geq 0. \]
\end{prop}
\begin{proof}
The proposition is certainly true when $\kappa_3 \geq \kappa_4 \geq 0$. We consider the other cases.

Since $\kappa_1 \geq \kappa_2 > 0$, we know that
\begin{equation} \label{eq5-5}
\kappa_2 \kappa_3 + \kappa_2 \kappa_4 + \kappa_1^2 \geq \kappa_2 \kappa_3 + \kappa_2 \kappa_4 + \kappa_1 \kappa_2 = \kappa_2 (\kappa_3 + \kappa_4 + \kappa_1) > 0.
\end{equation}
Then the case when $\kappa_4 \leq \kappa_3 \leq 0$ is proved.

For the case when $\kappa_3 \geq 0$ while $\kappa_4 \leq 0$, by the fact
\[ \kappa_2 + \kappa_3 > - \kappa_4 \geq 0, \]
we arrive at
\[ 2 \kappa_1^2 \geq \kappa_2 \kappa_3 + \kappa_3^2 \geq - \kappa_4 \kappa_3.  \]
Combined with \eqref{eq5-5} we can prove the proposition.
\end{proof}

\begin{prop} \label{Prop2}
\[ T_1 \geq 0. \]
\end{prop}
\begin{proof}

For the case when $\kappa_3 \geq \kappa_4 \geq 0$, by Proposition \ref{Prop0}, we immediately arrive at
\[ T_1 \geq 107 \kappa_1^6 \]
for $\kappa_1$ sufficiently large.

For the case when $\kappa_4 \leq \kappa_3 \leq 0$,
by Proposition \ref{Prop1} and Proposition \ref{Prop0}, for $\kappa_1$ sufficiently large,
\[ T_1 \geq 107 \kappa_1^6  - 72  k_1^6  + 108 \kappa_1^3 \sigma_3 + 12 \kappa_1^2 \sigma_2^2 - 36 \kappa_1^3 \sigma_3 \geq 0.   \]

The remaining case is when $\kappa_3 \geq 0$ while $\kappa_4 \leq 0$.
By Mathematica, we find that
\[ \frac{\partial^5 T_1}{\partial \kappa_1 \partial \kappa_4^4} = 72864 \kappa_1 + 68112 \kappa_2 + 68112 \kappa_3 + 46560 \kappa_4 > 0. \]
Then
\[ \begin{aligned}
 \frac{\partial^4 T_1}{\partial \kappa_1 \partial \kappa_4^3}
= & 43920 \kappa_1^2 + 84276 \kappa_1 \kappa_2 + 39636 \kappa_2^2 + 84276 \kappa_1 \kappa_3 + 81408 \kappa_2 \kappa_3 \\
& + 39636 \kappa_3^2 + 72864 \kappa_1 \kappa_4 + 68112 \kappa_2 \kappa_4
+ 68112 \kappa_3 \kappa_4 + 23280 \kappa_4^2 \\
\geq & 43920 \kappa_1^2 + 84276 \kappa_1 \kappa_2 + 39636 \kappa_2^2 + 84276 \kappa_1 \kappa_3 + 81408 \kappa_2 \kappa_3 \\
& + 39636 \kappa_3^2 - 72864 \kappa_1 (\kappa_2 + \kappa_3) - 68112 \kappa_2 (\kappa_2 + \kappa_3) \\
& - 68112 \kappa_3 (\kappa_2 + \kappa_3) + 23280 (\kappa_2 + \kappa_3)^2 \\
= & 12 (3660 \kappa_1^2 + 951 \kappa_1 \kappa_2 - 433 \kappa_2^2 + 951 \kappa_1 \kappa_3 - 688 \kappa_2 \kappa_3 - 433 \kappa_3^2) \geq 0.
\end{aligned} \]
Similarly, we can prove that
\[  \frac{\partial^3 T_1}{\partial \kappa_1 \partial \kappa_4^2} \geq 0, \qquad  \frac{\partial^2 T_1}{\partial \kappa_1 \partial \kappa_4} \geq 0. \]
Hence
\[ \begin{aligned}
\frac{\partial T_1}{\partial \kappa_1} & \geq 12 \Big( 54 \kappa_1^5 - 8 \kappa_1^3 \kappa_2^2 + 2 \kappa_1 \kappa_2^4 - 8 \kappa_1^3 \kappa_2 \kappa_3 \\
& - 27 \kappa_1^2 \kappa_2^2 \kappa_3 + 4 \kappa_1 \kappa_2^3 \kappa_3 + \kappa_2^4 \kappa_3 - 8 \kappa_1^3 \kappa_3^2 - 27 \kappa_1^2 \kappa_2 \kappa_3^2 + 6 \kappa_1 \kappa_2^2 \kappa_3^2 \\
& + 2 \kappa_2^3 \kappa_3^2 + 4 \kappa_1 \kappa_2 \kappa_3^3 + 2 \kappa_2^2 \kappa_3^3 + 2 \kappa_1 \kappa_3^4 + \kappa_2 \kappa_3^4 \Big) := U_1.
\end{aligned} \]
Since
\[ \begin{aligned}
\frac{\partial U_1}{\partial \kappa_1} = & 12 \Big( 270 \kappa_1^4 - 24 \kappa_1^2 \kappa_2^2 + 2 \kappa_2^4 - 24 \kappa_1^2 \kappa_2 \kappa_3 - 54 \kappa_1 \kappa_2^2 \kappa_3 + 4 \kappa_2^3 \kappa_3 \\
& - 24 \kappa_1^2 \kappa_3^2 - 54 \kappa_1 \kappa_2 \kappa_3^2 +
   6 \kappa_2^2 \kappa_3^2 + 4 \kappa_2 \kappa_3^3 + 2 \kappa_3^4 \Big) \geq 0,
\end{aligned} \]

\[ \begin{aligned}
\frac{\partial T_1}{\partial \kappa_1} \geq & U_1 \geq  12 \Big( 54 \kappa_2^5 - 8 \kappa_2^3 \kappa_2^2 + 2 \kappa_2 \kappa_2^4 - 8 \kappa_2^3 \kappa_2 \kappa_3 \\
& - 27 \kappa_2^2 \kappa_2^2 \kappa_3 + 4 \kappa_2 \kappa_2^3 \kappa_3 + \kappa_2^4 \kappa_3 - 8 \kappa_2^3 \kappa_3^2 - 27 \kappa_2^2 \kappa_2 \kappa_3^2 + 6 \kappa_2 \kappa_2^2 \kappa_3^2 \\
& + 2 \kappa_2^3 \kappa_3^2 + 4 \kappa_2 \kappa_2 \kappa_3^3 + 2 \kappa_2^2 \kappa_3^3 + 2 \kappa_2 \kappa_3^4 + \kappa_2 \kappa_3^4 \Big) \\
= & 36 \kappa_2 (\kappa_2 - \kappa_3) (2 \kappa_2 + \kappa_3) (8 \kappa_2^2 - \kappa_2 \kappa_3 - \kappa_3^2) \geq 0.
\end{aligned} \]
It follows that
\begin{equation*}
\begin{aligned}
& T_1 \geq T_1 \Big|_{\kappa_1 = \kappa_2} =  24 \kappa_2^2 (\kappa_2 - \kappa_3)^2 (2 \kappa_2 + \kappa_3)^2 \\
& + \sigma_1 \Bigg( 2 \kappa_2 (\kappa_2 - \kappa_3) (2 \kappa_2 + \kappa_3) (103 \kappa_2^2 + 28 \kappa_2 \kappa_3 + 4 \kappa_3^2) \\
 & + \sigma_1 \bigg( 1304 \kappa_2^4 - 171 \kappa_2^3 \kappa_3 - 285 \kappa_2^2 \kappa_3^2 + 88 \kappa_2 \kappa_3^3 + 36 \kappa_3^4  \\
 & + \sigma_1 \Big( 2453 \kappa_2^3 - 349 \kappa_2^2 \kappa_3 - 790 \kappa_2 \kappa_3^2 - 72 \kappa_3^3 + \sigma_1 \big( 2181 \kappa_2^2 + 439 \kappa_2 \kappa_3 - 235 \kappa_3^2 \\
 & + \sigma_1 (674 \kappa_2 + 286 \kappa_3 + 15 \kappa_4) \big) \Big) \bigg) \Bigg)  \geq 0.
\end{aligned}
\end{equation*}

\end{proof}

\vspace{4mm}

\section{For $i = 2$}

\vspace{4mm}

By Mathematica, we find
\begin{footnotesize}
\begin{equation*}
\begin{aligned}
Q_2 = & - \sum\limits_{j l} f_{jl}  h_{jj2} h_{ll2} +  \frac{2 \sigma}{4} \frac{1}{H - \kappa_2}  \frac{1}{H - \kappa_1} t_1^2 + \frac{2 \sigma}{4} \frac{1}{H - \kappa_2} \frac{1}{H - \kappa_3} t_3^2 \\
& + \frac{2 \sigma}{4} \frac{1}{H - \kappa_2} \frac{1}{H - \kappa_4} t_4^2 \\
= &  \frac{\sigma}{4} \frac{\beta^2}{(\nu^{4 + 1})^2} \frac{u_2^2}{u^2} (\nu^{4 + 1} - \kappa_2)^2 H^2 \frac{1}{(\kappa_1 + \kappa_2 + \kappa_3) (\kappa_1 + \kappa_2 + \kappa_4) (\kappa_2 + \kappa_3 +
     \kappa_4)} \\
& \frac{2 \left( \substack{ 15 \kappa_1^5 + 72 \kappa_2^5 + 15 \kappa_3^5 + 196 \kappa_3^4 \kappa_4 + 509 \kappa_3^3 \kappa_4^2 +
     509 \kappa_3^2 \kappa_4^3 + 196 \kappa_3 \kappa_4^4 + 15 \kappa_4^5 + 348 \kappa_2^4 (\kappa_3 + \kappa_4) \\
     + 2 \kappa_2^3 (325 \kappa_3^2 + 668 \kappa_3 \kappa_4 + 325 \kappa_4^2) +
     3 \kappa_2^2 (189 \kappa_3^3 + 617 \kappa_3^2 \kappa_4 + 617 \kappa_3 \kappa_4^2 + 189 \kappa_4^3) \\
     + \kappa_2 (208 \kappa_3^4 + 1063 \kappa_3^3 \kappa_4 + 1722 \kappa_3^2 \kappa_4^2 + 1063 \kappa_3 \kappa_4^3 +
        208 \kappa_4^4) \\
         + 4 \kappa_1^4 (52 \kappa_2 + 49 (\kappa_3 + \kappa_4)) +
     \kappa_1^3 (567 \kappa_2^2 + 509 \kappa_3^2 + 1029 \kappa_3 \kappa_4 + 509 \kappa_4^2 +
        1063 \kappa_2 (\kappa_3 + \kappa_4)) \\
        + \kappa_1^2 (650 \kappa_2^3 + 509 \kappa_3^3 + 1662 \kappa_3^2 \kappa_4 + 1662 \kappa_3 \kappa_4^2 + 509 \kappa_4^3 + 1851 \kappa_2^2 (\kappa_3 + \kappa_4) +
        6 \kappa_2 (287 \kappa_3^2 + 598 \kappa_3 \kappa_4 + 287 \kappa_4^2)) \\
        + \kappa_1 (348 \kappa_2^4 + 196 \kappa_3^4 + 1029 \kappa_3^3 \kappa_4 + 1662 \kappa_3^2 \kappa_4^2 + 1029 \kappa_3 \kappa_4^3 + 196 \kappa_4^4 + 1336 \kappa_2^3 (\kappa_3 + \kappa_4) \\
        + 3 \kappa_2^2 (617 \kappa_3^2 + 1296 \kappa_3 \kappa_4 + 617 \kappa_4^2)
        + \kappa_2 (1063 \kappa_3^3 + 3588 \kappa_3^2 \kappa_4 + 3588 \kappa_3 \kappa_4^2 +
           1063 \kappa_4^3)) } \right)}{ \left( \substack{ 15 \kappa_1^4 + 12 \kappa_2^4 + 20 \kappa_2^3 (\kappa_3 + \kappa_4) + \kappa_2^2 (-5 \kappa_3^2 + 2 \kappa_3 \kappa_4 - 5 \kappa_4^2) + (\kappa_3 + \kappa_4)^2 (15 \kappa_3^2 - 14 \kappa_3 \kappa_4 + 15 \kappa_4^2) \\
           - 2 \kappa_2 (\kappa_3^3 + 11 \kappa_3^2 \kappa_4 + 11 \kappa_3 \kappa_4^2 + \kappa_4^3) -
     2 \kappa_1^3 (\kappa_2 - 8 (\kappa_3 + \kappa_4)) +
     \kappa_1^2 (-5 \kappa_2^2 - 22 \kappa_2 (\kappa_3 + \kappa_4) + 2 (\kappa_3^2 + \kappa_4^2)) \\
      + 2 \kappa_1 (10 \kappa_2^3 + \kappa_2^2 (\kappa_3 + \kappa_4) -
        \kappa_2 (11 \kappa_3^2 + 36 \kappa_3 \kappa_4 + 11 \kappa_4^2) + 8 (\kappa_3^3 + \kappa_4^3)) } \right) }.
\end{aligned}
\end{equation*}
\end{footnotesize}

\vspace{2mm}

\subsection{The case when $0 < \kappa_2 \leq \nu^{n + 1}$.}~

\vspace{2mm}

In this case, we have
\[   \kappa_4 \leq \kappa_3 \leq \kappa_2 \leq \nu^{n + 1}, \]
and so
\[ \kappa_4 > - \kappa_2 - \kappa_3 \geq - 2 \nu^{n + 1}. \]
Hence $\kappa_2$, $\kappa_3$ and $\kappa_4$ are uniformly bounded.

In \eqref{eq5-2}, we want to guarantee the nonnegativity of the following $i = 2$ terms
\begin{equation*}
\begin{aligned}
& ( \beta - 1 ) H \bigg( f_2 + f_2 \kappa_2^2 \bigg) + Q_2 + \frac{2 \beta}{\nu^{n + 1}} H  f_2 \frac{u_2^2}{u^2} \big( \kappa_2 - \nu^{n+1} \big) \\
& - \frac{\beta (\beta - 1)}{(\nu^{n + 1})^2} H  f_2  \frac{u_2^2}{u^2} (\nu^{n + 1} - \kappa_2)^2.
\end{aligned}
\end{equation*}
We note that
\[ \begin{aligned}
\frac{Q_2}{H f_2} \geq (2 - \delta) \frac{\beta^2}{(\nu^{n + 1})^2} \frac{u_2^2}{u^2} (\kappa_2 - \nu^{n + 1})^2,
\end{aligned} \]
where $\delta$ is a small positive constant which is to be chosen later.

It suffices to guarantee the positivity of
\begin{equation*}
\begin{aligned}
& \bigg( 1 + \frac{1}{\beta} - \delta \bigg) \frac{\beta^2}{(\nu^{n + 1})^2} (\kappa_2 - \nu^{n + 1})^2 + \frac{2 \beta}{\nu^{n + 1}} \big( \kappa_2 - \nu^{n+1} \big) + ( \beta - 1 ).
\end{aligned}
\end{equation*}
For this, we require that
\begin{equation} \label{Condition1}
1 - (\beta - 1)  \bigg( 1 + \frac{1}{\beta} - \delta \bigg) < 0.
\end{equation}
We first require that
\begin{equation*}
1 - (\beta - 1)  \bigg( 1 + \frac{1}{\beta} \bigg) < 0,
\end{equation*}
which implies that
\begin{equation} \label{Condition1-1}
\beta > \frac{1 + \sqrt{5}}{2}.
\end{equation}

\vspace{2mm}

\subsection{The case when $\kappa_2 \geq \nu^{n + 1}$.}~

\vspace{2mm}

For any $\beta > 1$,

\begin{footnotesize}
\begin{equation*}
\begin{aligned}
& - \sum\limits_{j l} f_{jl}  h_{jj2} h_{ll2} +   \frac{2 \sigma}{4} \frac{1}{H - \kappa_2} \frac{1}{H - \kappa_1} t_1^2 + \frac{2 \sigma}{4} \frac{1}{H - \kappa_2} \frac{1}{H - \kappa_3} t_3^2 \\
& + \frac{2 \sigma}{4} \frac{1}{H - \kappa_2} \frac{1}{H - \kappa_4} t_4^2
- \frac{\beta^2}{(\nu^{4 + 1})^2} H  f_2  \frac{u_2^2}{u^2} (\kappa_2 - \nu^{4 + 1})^2 \\
= &  \frac{\sigma}{4} \frac{\beta^2}{(\nu^{4 + 1})^2} \frac{u_2^2}{u^2} (\kappa_2 - \nu^{4 + 1})^2 H \frac{1}{(\kappa_1 + \kappa_2 + \kappa_3) (\kappa_1 + \kappa_2 +
     \kappa_4) (\kappa_2 + \kappa_3 + \kappa_4)} \\
& \cdot \frac{T_2}{\left(  \substack{ 15 \kappa_1^4 -
   2 \kappa_1^3 \big( \kappa_2 - 8 (\kappa_3 + \kappa_4) \big) +
   \kappa_1^2 \big( -5 \kappa_2^2 - 22 \kappa_2 (\kappa_3 + \kappa_4) + 2 (\kappa_3^2 + \kappa_4^2) \big) \\ + 2 \kappa_1 \big( 10 \kappa_2^3 + \kappa_2^2 (\kappa_3 + \kappa_4) -
      \kappa_2 (11 \kappa_3^2 + 36 \kappa_3 \kappa_4 + 11 \kappa_4^2) + 8 (\kappa_3^3 + \kappa_4^3) \big)
   + 12 \kappa_2^4 + 20 \kappa_2^3 (\kappa_3 + \kappa_4) \\
   +  \kappa_2^2 (-5 \kappa_3^2 + 2 \kappa_3 \kappa_4 - 5 \kappa_4^2) + (\kappa_3 + \kappa_4)^2 (15 \kappa_3^2 - 14 \kappa_3 \kappa_4 + 15 \kappa_4^2) -
   2 \kappa_2 (\kappa_3^3 + 11 \kappa_3^2 \kappa_4 + 11 \kappa_3 \kappa_4^2 + \kappa_4^3) }  \right)},
\end{aligned}
\end{equation*}
\end{footnotesize}
where
\begin{equation*}
\begin{aligned}
T_2 := & 15 \kappa_1^6 + \kappa_1^5 \big( 388 \kappa_2 + 361 (\kappa_3 + \kappa_4) \big) \\
& + \kappa_1^4 \big( 1518 \kappa_2^2 + 1345 \kappa_3^2 + 2701 \kappa_3 \kappa_4 + 1345 \kappa_4^2 +
      2838 \kappa_2 (\kappa_3 + \kappa_4) \big) \\
& + \kappa_1^3 \Big( 2440 \kappa_2^3 + 7023 \kappa_2^2 (\kappa_3 + \kappa_4) +
      2 \kappa_2 (3303 \kappa_3^2 + 6784 \kappa_3 \kappa_4 + 3303 \kappa_4^2) \\
      & +   6 (333 \kappa_3^3 + 1055 \kappa_3^2 \kappa_4 + 1055 \kappa_3 \kappa_4^2 + 333 \kappa_4^3) \Big) \\
& + \kappa_1^2 \Big( 1919 \kappa_2^4 + 1345 \kappa_3^4 + 6330 \kappa_3^3 \kappa_4 + 9966 \kappa_3^2 \kappa_4^2 + 6330 \kappa_3 \kappa_4^3 \\
& + 1345 \kappa_4^4 + 7672 \kappa_2^3 (\kappa_3 + \kappa_4)  + 3 \kappa_2^2 (3674 \kappa_3^2 + 7607 \kappa_3 \kappa_4 + 3674 \kappa_4^2) \\
& +  6 \kappa_2 (1101 \kappa_3^3 + 3584 \kappa_3^2 \kappa_4 + 3584 \kappa_3 \kappa_4^2 + 1101 \kappa_4^3) \Big) \\
& + \kappa_1 \Big( 732 \kappa_2^5 + 361 \kappa_3^5 + 2701 \kappa_3^4 \kappa_4 + 6330 \kappa_3^3 \kappa_4^2 +
      6330 \kappa_3^2 \kappa_4^3 + 2701 \kappa_3 \kappa_4^4 \\
 & + 361 \kappa_4^5 + 3862 \kappa_2^4 (\kappa_3 + \kappa_4)
 + 4 \kappa_2^3 (1918 \kappa_3^2 + 3951 \kappa_3 \kappa_4 + 1918 \kappa_4^2) \\
  &    +  3 \kappa_2^2 (2341 \kappa_3^3 + 7607 \kappa_3^2 \kappa_4 + 7607 \kappa_3 \kappa_4^2 + 2341 \kappa_4^3) \\
& + 2 \kappa_2 (1419 \kappa_3^4 + 6784 \kappa_3^3 \kappa_4 + 10752 \kappa_3^2 \kappa_4^2 +
         6784 \kappa_3 \kappa_4^3 + 1419 \kappa_4^4) \Big) \\
& + 108 \kappa_2^6 + 732 \kappa_2^5 (\kappa_3 + \kappa_4) +
   \kappa_2^4 (1919 \kappa_3^2 + 3862 \kappa_3 \kappa_4 + 1919 \kappa_4^2) \\
& +  8 \kappa_2^3 (305 \kappa_3^3 + 959 \kappa_3^2 \kappa_4 + 959 \kappa_3 \kappa_4^2 + 305 \kappa_4^3) \\
&   + (\kappa_3 + \kappa_4)^2 (15 \kappa_3^4 + 331 \kappa_3^3 \kappa_4 + 668 \kappa_3^2 \kappa_4^2 + 331 \kappa_3 \kappa_4^3 + 15 \kappa_4^4) \\
& + 3 \kappa_2^2 (506 \kappa_3^4 + 2341 \kappa_3^3 \kappa_4 + 3674 \kappa_3^2 \kappa_4^2 + 2341 \kappa_3 \kappa_4^3 + 506 \kappa_4^4) \\
&   + \kappa_2 (388 \kappa_3^5 + 2838 \kappa_3^4 \kappa_4 + 6606 \kappa_3^3 \kappa_4^2 + 6606 \kappa_3^2 \kappa_4^3 + 2838 \kappa_3 \kappa_4^4 + 388 \kappa_4^5).
\end{aligned}
\end{equation*}

\begin{prop} \label{T2}
\[ T_2 \geq 0. \]
\end{prop}
\begin{proof}
For the case when $\kappa_3 \geq \kappa_4 \geq 0$, we can immediately see that $T_2 \geq 0$.

For the case when $\kappa_4 \leq \kappa_3 \leq 0$, we can write
$T_2$ as
\begin{equation*}
\begin{aligned}
T_2 = & (\kappa_1 + \kappa_2) (\kappa_1 - \kappa_3) (\kappa_1 + \kappa_2 + \kappa_3) \Big( 15 \kappa_1^3 - 3 \kappa_1^2 \kappa_2 \\
& + 16 \kappa_1 \kappa_2^2 - 4 \kappa_2^3 + 4 \kappa_1 \kappa_2 \kappa_3 - 16 \kappa_2^2 \kappa_3 + 4 \kappa_1 \kappa_3^2 - 16 \kappa_2 \kappa_3^2 \Big) \\
& + \sigma_1 P_2,
\end{aligned}
\end{equation*}
where $P_2$ is a homogeneous polynomial of degree $5$.
By Proposition \ref{Prop0}, we know that for $\kappa_1$ sufficiently large,
\begin{equation*}
\begin{aligned}
T_2 \geq & \kappa_1  \kappa_1  \kappa_1 \Big( 15 \kappa_1^3 - 3 \kappa_1^2 \kappa_2 + 15 \kappa_1 \kappa_2^2 - 4 \kappa_2^3 \\
& - 16 \kappa_2 \kappa_3 (\kappa_3 + \kappa_2) + \kappa_1 (\kappa_2 + 2 \kappa_3)^2 \Big)  - \kappa_1^6 \geq 0.
\end{aligned}
\end{equation*}

For the case when $\kappa_3 \geq 0$ while $\kappa_4 \leq 0$,
we prove the nonnegativity of $T_2$ similarly as $T_1$. By Mathematica,
\begin{equation*}
\frac{\partial^5 T_2}{\partial \kappa_1 \partial \kappa_4^4} = 24 (2690 \kappa_1 + 2838 \kappa_2 + 2701 \kappa_3 + 1805 \kappa_4) \geq 0.
\end{equation*}
Hence
\begin{equation*}
\begin{aligned}
& \frac{\partial^4 T_2}{\partial \kappa_1 \partial \kappa_4^3} \geq \frac{\partial^4 T_2}{\partial \kappa_1 \partial \kappa_4^3} \Big|_{\kappa_4 = - \kappa_2 - \kappa_3} \\
= & 6 \Big( 5994 \kappa_1^2 - 719 \kappa_2^2 - 1368 \kappa_2 \kappa_3 - 864 \kappa_3^2 +
   4 \kappa_1 (613 \kappa_2 + 475 \kappa_3) \Big) \geq 0.
\end{aligned}
\end{equation*}
Then we have
\begin{equation*}
\begin{aligned}
& \frac{\partial^3 T_2}{\partial \kappa_1 \partial \kappa_4^2} \geq \frac{\partial^3 T_2}{\partial \kappa_1 \partial \kappa_4^2} \Big|_{\kappa_4 = - \kappa_2 - \kappa_3} \geq 0.
\end{aligned}
\end{equation*}
It follows that
\begin{equation*}
\begin{aligned}
& \frac{\partial^2 T_2}{\partial \kappa_1 \partial \kappa_4} \geq \frac{\partial^2 T_2}{\partial \kappa_1 \partial \kappa_4} \Big|_{\kappa_4 = - \kappa_2 - \kappa_3} \\
= & 1805 \kappa_1^4 + 592 \kappa_1^3 \kappa_2 - 585 \kappa_1^2 \kappa_2^2 + 132 \kappa_1 \kappa_2^3 + 40 \kappa_2^4 + 44 \kappa_1^3 \kappa_3 \\
& - 948 \kappa_1^2 \kappa_2 \kappa_3  + 510 \kappa_1 \kappa_2^2 \kappa_3 + 20 \kappa_2^3 \kappa_3 -
 1008 \kappa_1^2 \kappa_3^2 \\
& + 444 \kappa_1 \kappa_2 \kappa_3^2 + 16 \kappa_1 \kappa_3^3 + 40 \kappa_2 \kappa_3^3 + 32 \kappa_3^4 := U_2.
\end{aligned}
\end{equation*}
Since
\[ \begin{aligned}
\frac{\partial U_2}{\partial \kappa_1} = & 7220 \kappa_1^3 + 1776 \kappa_1^2 \kappa_2 - 1170 \kappa_1 \kappa_2^2 + 132 \kappa_2^3 + 132 \kappa_1^2 \kappa_3 \\
& - 1896 \kappa_1 \kappa_2 \kappa_3 + 510 \kappa_2^2 \kappa_3 - 2016 \kappa_1 \kappa_3^2 + 444 \kappa_2 \kappa_3^2 + 16 \kappa_3^3 \geq 0,
\end{aligned} \]
hence
\[ \begin{aligned}
\frac{\partial^2 T_2}{\partial \kappa_1 \partial \kappa_4} \geq & U_2 \geq U_2 \Big|_{\kappa_1 = \kappa_2} \\
= & 2 \Big( 992 \kappa_2^4 - 187 \kappa_2^3 \kappa_3 - 282 \kappa_2^2 \kappa_3^2 + 28 \kappa_2 \kappa_3^3 + 16 \kappa_3^4 \Big) \geq 0,
\end{aligned} \]
and so
\[ \begin{aligned}
& \frac{\partial T_2}{\partial \kappa_1} \geq \frac{\partial T_2}{\partial \kappa_1}\Big|_{\kappa_4 = - \kappa_2 - \kappa_3} \\
= & 90 \kappa_1^5 + 135 \kappa_1^4 \kappa_2 + 100 \kappa_1^3 \kappa_2^2 + 75 \kappa_1^2 \kappa_2^3 + 16 \kappa_1 \kappa_2^4 - 4 \kappa_2^5 - 44 \kappa_1^3 \kappa_2 \kappa_3 \\
& - 60 \kappa_1^2 \kappa_2^2 \kappa_3 - 82 \kappa_1 \kappa_2^3 \kappa_3 -
 28 \kappa_2^4 \kappa_3 - 44 \kappa_1^3 \kappa_3^2 - 60 \kappa_1^2 \kappa_2 \kappa_3^2 \\
& - 90 \kappa_1 \kappa_2^2 \kappa_3^2 -
 16 \kappa_2^3 \kappa_3^2 - 16 \kappa_1 \kappa_2 \kappa_3^3 + 24 \kappa_2^2 \kappa_3^3 - 8 \kappa_1 \kappa_3^4 + 12 \kappa_2 \kappa_3^4 := V_2.
\end{aligned} \]
Since
\[ \begin{aligned}
\frac{\partial V_2}{\partial \kappa_1} = & 450 \kappa_1^4 + 540 \kappa_1^3 \kappa_2 + 300 \kappa_1^2 \kappa_2^2 + 150 \kappa_1 \kappa_2^3 + 16 \kappa_2^4 \\
& - 132 \kappa_1^2 \kappa_2 \kappa_3 - 120 \kappa_1 \kappa_2^2 \kappa_3 - 82 \kappa_2^3 \kappa_3 - 132 \kappa_1^2 \kappa_3^2 \\
& - 120 \kappa_1 \kappa_2 \kappa_3^2 - 90 \kappa_2^2 \kappa_3^2 - 16 \kappa_2 \kappa_3^3 - 8 \kappa_3^4 \geq 0,
\end{aligned} \]
hence
\[ \begin{aligned}
V_2 \geq & V_2 \Big|_{\kappa_1 = \kappa_2} \\
= & 2 \Big( 206 \kappa_2^5 - 107 \kappa_2^4 \kappa_3 - 105 \kappa_2^3 \kappa_3^2 + 4 \kappa_2^2 \kappa_3^3 + 2 \kappa_2 \kappa_3^4 \Big) := W_2.
\end{aligned} \]
Since
\[ \frac{\partial W_2}{\partial \kappa_2} = 2 \Big( 1030 \kappa_2^4 - 428 \kappa_2^3 \kappa_3 - 315 \kappa_2^2 \kappa_3^2 + 8 \kappa_2 \kappa_3^3 + 2 \kappa_3^4 \Big) \geq 0, \]
we know that
\[ \begin{aligned}
& W_2 \geq W_2 \Big|_{\kappa_2 = \kappa_3}
=  0.
\end{aligned} \]
Hence
\[ \frac{\partial T_2}{\partial \kappa_1} \geq 0. \]
Therefore we have
\begin{equation*}
\begin{aligned}
T_2 \geq & T_2 \Big|_{\kappa_1 = \kappa_2} = 24 \kappa_2^2 (\kappa_2 - \kappa_3)^2 (2 \kappa_2 + \kappa_3)^2 \\
& + \sigma_1 \Bigg( 2 \kappa_2 (\kappa_2 - \kappa_3) (2 \kappa_2 + \kappa_3) (103 \kappa_2^2 + 28 \kappa_2 \kappa_3 + 4 \kappa_3^2) \\
& + \sigma_1 \bigg( 1304 \kappa_2^4 - 171 \kappa_2^3 \kappa_3 - 285 \kappa_2^2 \kappa_3^2 + 88 \kappa_2 \kappa_3^3 + 36 \kappa_3^4 \\
& + \sigma_1 \Big( 2453 \kappa_2^3 - 349 \kappa_2^2 \kappa_3 - 790 \kappa_2 \kappa_3^2 - 72 \kappa_3^3 \\
& + \sigma_1 \big( 2181 \kappa_2^2 + 439 \kappa_2 \kappa_3 - 235 \kappa_3^2 + \sigma_1 (674 \kappa_2 + 286 \kappa_3 + 15 \kappa_4) \big) \Big) \bigg) \Bigg) \geq 0.
\end{aligned}
\end{equation*}
Here
\[ \sigma_1 = \kappa_2 + \kappa_3 + \kappa_4. \]
\end{proof}

\vspace{4mm}

\section{For $i = 3$}

\vspace{4mm}

In this section, we divide our discussion into two cases. That is, the case when $\kappa_3 \geq - \varphi \kappa_1$ and the case when $\kappa_3 \leq - \varphi \kappa_1$. Here $0 < \varphi < 1$ is a constant to be determined later.

\vspace{2mm}

\subsection{The case when $\kappa_3 \geq - \varphi \kappa_1$.}~

\vspace{2mm}

Recall that
\begin{equation*}
\begin{aligned}
Q_3 = & - \sum\limits_{j l} f_{jl}  h_{jj3} h_{ll3} +  \frac{2 \sigma}{4} \frac{1}{H - \kappa_3}  \frac{1}{H - \kappa_1} t_1^2 + \frac{2 \sigma}{4} \frac{1}{H - \kappa_3} \frac{1}{H - \kappa_2} t_2^2 \\
& + \frac{2 \sigma}{4} \frac{1}{H - \kappa_3} \frac{1}{H - \kappa_4} t_4^2 .
\end{aligned}
\end{equation*}

Since
\[  \frac{2 \beta}{\nu^{4 + 1}} H f_3 \frac{u_3^2}{u^2} \big( \kappa_3 - \nu^{4 + 1} \big) \geq - \frac{\beta^2}{(\beta - 1)(\nu^{4 + 1})^2} H f_3 \frac{u_3^2}{u^2} (\kappa_3 - \nu^{4 + 1})^2 - (\beta - 1) f_3 \frac{u_3^2}{u^2} H,   \]
we want to guarantee the nonnegativity of the following quantity:
\begin{equation*}
\begin{aligned}
& Q_3 - \frac{\beta^2}{(\beta - 1)(\nu^{4 + 1})^2} H f_3 \frac{u_3^2}{u^2} (\kappa_3 - \nu^{4 + 1})^2 - \frac{\beta (\beta - 1)}{(\nu^{4 + 1})^2} H f_3 \frac{u_3^2}{u^2} (\kappa_3 - \nu^{4 + 1})^2 \\
= &  \frac{\sigma}{4} \frac{\beta}{(\beta - 1)(\nu^{4 + 1})^2} \frac{u_3^2}{u^2} (\kappa_3 - \nu^{4 + 1})^2 H \frac{T_3}{ (\kappa_1 + \kappa_2 + \kappa_3) (\kappa_1 + \kappa_3 + \kappa_4) (\kappa_2 + \kappa_3 + \kappa_4) B_3},
\end{aligned}
\end{equation*}
where
\begin{equation*}
\begin{aligned}
B_3 := & 15 \kappa_1^4 + 16 \kappa_1^3 \kappa_2 + 2 \kappa_1^2 \kappa_2^2 + 16 \kappa_1 \kappa_2^3 +
   15 \kappa_2^4 - 2 \kappa_1^3 \kappa_3 - 22 \kappa_1^2 \kappa_2 \kappa_3 \\
   & - 22 \kappa_1 \kappa_2^2 \kappa_3 - 2 \kappa_2^3 \kappa_3 -
   5 \kappa_1^2 \kappa_3^2 + 2 \kappa_1 \kappa_2 \kappa_3^2 - 5 \kappa_2^2 \kappa_3^2 + 20 \kappa_1 \kappa_3^3 \\
   & + 20 \kappa_2 \kappa_3^3 + 12 \kappa_3^4 + 16 \kappa_1^3 \kappa_4 + 16 \kappa_2^3 \kappa_4 - 22 \kappa_1^2 \kappa_3 \kappa_4 - 72 \kappa_1 \kappa_2 \kappa_3 \kappa_4 \\
 &  - 22 \kappa_2^2 \kappa_3 \kappa_4 + 2 \kappa_1 \kappa_3^2 \kappa_4 + 2 \kappa_2 \kappa_3^2 \kappa_4 + 20 \kappa_3^3 \kappa_4 + 2 \kappa_1^2 \kappa_4^2 + 2 \kappa_2^2 \kappa_4^2 - 22 \kappa_1 \kappa_3 \kappa_4^2 \\ & - 22 \kappa_2 \kappa_3 \kappa_4^2 - 5 \kappa_3^2 \kappa_4^2 + 16 \kappa_1 \kappa_4^3 + 16 \kappa_2 \kappa_4^3 - 2 \kappa_3 \kappa_4^3 + 15 \kappa_4^4,
\end{aligned}
\end{equation*}
\begin{equation} \label{eq6-2}
\begin{aligned}
& T_3 := (\kappa_1 - \kappa_2) (\kappa_1 + \kappa_2 + \kappa_3) \Big( (\beta^2 - \beta - 1) \kappa_1 + (\beta^2 - \beta + 1) \kappa_3 \Big) \\
& \cdot \Big( 15 \kappa_1^3 - 3 \kappa_1^2 \kappa_3  - 4 \kappa_3 (2 \kappa_2 + \kappa_3)^2 + 4 \kappa_1 (\kappa_2^2 + \kappa_2 \kappa_3 + 4 \kappa_3^2) \Big) + (\kappa_2 + \kappa_3 + \kappa_4) M_3,
\end{aligned}
\end{equation}
and $M_3$ is a homogeneous polynomial of degree $5$.

We differentiate $M_3$ five times with respect to $\kappa_4$ by Mathematica. We find that
\[ \frac{\partial^5 M_3}{\partial \kappa_4^5} = 1800 ( \beta^2 - \beta - 1) > 0 \]
by Condition \eqref{Condition1-1}. Hence
\[ \begin{aligned}
& \frac{\partial^4 M_3}{\partial \kappa_4^4} \geq \frac{\partial^4 M_3}{\partial \kappa_4^4} \Big|_{\kappa_4 = - \kappa_2 - \kappa_3} \\
= & 24 \Big( (-61 - 361 \beta + 361 \beta^2) \kappa_1 + (29 - 271 \beta + 271 \beta^2) \kappa_2 +
   2 (16 - 149 \beta + 149 \beta^2) \kappa_3 \Big) \\
\geq & 0.
\end{aligned} \]
It follows that
\[ \begin{aligned}
& \frac{\partial^3 M_3}{\partial \kappa_4^3} \geq \frac{\partial^3 M_3}{\partial \kappa_4^3} \Big|_{\kappa_4 = - \kappa_2 - \kappa_3} \\
= & 6 \Big( 5 (-13 - 269 \beta + 269 \beta^2) \kappa_1^2 + 4 (41 - 224 \beta + 224 \beta^2) \kappa_1 \kappa_2 - 5 (-3 - 47 \beta + 47 \beta^2) \kappa_2^2 \\
& + (209 - 1033 \beta + 1033 \beta^2) \kappa_1 \kappa_3 + (49 + 457 \beta - 457 \beta^2) \kappa_2 \kappa_3 + (33 + 197 \beta - 197 \beta^2) \kappa_3^2 \Big) \\
:= & N_3.
\end{aligned} \]
Since
\[ \begin{aligned}
& \frac{\partial N_3}{\partial \kappa_1}
=  6 \Big( 10 (-13 - 269 \beta + 269 \beta^2) \kappa_1 \\
& + 4 (41 - 224 \beta + 224 \beta^2) \kappa_2 + (209 - 1033 \beta + 1033 \beta^2) \kappa_3 \Big) \geq 0,
\end{aligned} \]
we know that
\[ \begin{aligned}
N_3 \geq & N_3 \Big|_{\kappa_1 = \kappa_2}
=  6 \Big( 2 (57 - 1003 \beta + 1003 \beta^2) \kappa_2^2 \\
& + 6 (43 - 96 \beta + 96 \beta^2) \kappa_2 \kappa_3 + (33 + 197 \beta - 197 \beta^2) \kappa_3^2 \Big) \geq 0.
\end{aligned} \]
Similarly, we can prove that
\[ \frac{\partial^2 M_3}{\partial \kappa_4^2} \geq 0, \qquad  \frac{\partial M_3}{\partial \kappa_4} \geq 0. \]
Therefore,
\[ M_3 \geq  M_3 \Big|_{\kappa_4 = - \kappa_2 - \kappa_3} := U_3. \]
Similar as above, by differentiating $U_3$ several times with respect to $\kappa_1$, we can find that
\[ \frac{\partial^2 U_3}{\partial \kappa_1^2} \geq 0. \]
Hence
\[ \begin{aligned}
& \frac{\partial U_3}{\partial \kappa_1} \geq  \frac{\partial U_3}{\partial \kappa_1} \Big|_{\kappa_1 = \kappa_2} \\
= & (-449 - 889 \beta + 889 \beta^2) \kappa_2^4 +
   2 (115 - 64 \beta + 64 \beta^2) \kappa_2^3 \kappa_3 +
   3 (11 + 25 \beta - 25 \beta^2) \kappa_2^2 \kappa_3^2 \\
 &  + 2 (61 - 76 \beta + 76 \beta^2) \kappa_2 \kappa_3^3 + 8 (8 - 5 \beta + 5 \beta^2) \kappa_3^4 := V_3.
\end{aligned} \]
We realize that
\[ \frac{\partial V_3}{\partial \kappa_2} \geq 0.  \]
Hence
\[ V_3 \geq V_3 \Big|_{\kappa_2 = - \kappa_3} = 2 (-352 - 287 \beta + 287 \beta^2) \kappa_3^4. \]

Now we impose the following condition:
\begin{equation} \label{Condition5}
\beta > \beta_0,
\end{equation}
where $\beta_0$ is the positive root of the quadratic polynomial
\[ -352 - 287 \beta + 287 \beta^2. \]
Mathematica shows that
\[ \beta_0 \approx 1.71511. \]

Then we know that
\begin{equation} \label{eq6-1}
\begin{aligned}
U_3 \geq  U_3 \Big|_{\kappa_1 = \kappa_2}
= & 2 (\kappa_2 - \kappa_3) (2 \kappa_2 + \kappa_3) \Big( (\beta^2 - \beta - 1) \kappa_2 + (\beta^2 - \beta + 1) \kappa_3 \Big) \\
& \cdot (19 \kappa_2^2 + 4 \kappa_2 \kappa_3 + 4 \kappa_3^2).
\end{aligned}
\end{equation}

In view of \eqref{eq6-2} and \eqref{eq6-1}, we divide our discussion into two subcases:

\begin{enumerate}
  \item The subcase when  $\kappa_2 \geq \phi \kappa_1$;
  \item The subcase when  $\kappa_2 \leq \phi \kappa_1$.
\end{enumerate}
Here $\phi \in (\varphi, 1)$ is a constant to be determined later.

{\bf Subcase (1).} In view of \eqref{eq6-1}, we require that
\[  (\beta^2 - \beta - 1) \phi - (\beta^2 - \beta + 1) \varphi > 0, \]
or equivalently,
\begin{equation} \label{Condition6}
 \frac{\phi}{\varphi} > \frac{\beta^2 - \beta + 1}{\beta^2 - \beta - 1}.
\end{equation}
Therefore, $M_3 \geq 0$. By \eqref{eq6-2}, we know that
\begin{equation} \label{eq6-4}
\begin{aligned}
& T_3 \geq (\kappa_1 - \kappa_2) (\kappa_1 + \kappa_2 + \kappa_3) \Big( (\beta^2 - \beta - 1) \kappa_1 + (\beta^2 - \beta + 1) \kappa_3 \Big) \\
& \cdot \Big( 15 \kappa_1^3 - 3 \kappa_1^2 \kappa_3  - 4 \kappa_3 (2 \kappa_2 + \kappa_3)^2 + 4 \kappa_1 (\kappa_2^2 + \kappa_2 \kappa_3 + 4 \kappa_3^2) \Big).
\end{aligned}
\end{equation}
Denote
\[ W_3 :=  15 \kappa_1^3 - 3 \kappa_1^2 \kappa_3  - 4 \kappa_3 (2 \kappa_2 + \kappa_3)^2 + 4 \kappa_1 (\kappa_2^2 + \kappa_2 \kappa_3 + 4 \kappa_3^2). \]
Since
\[ \frac{\partial W_3}{\partial \kappa_1} = 45 \kappa_1^2 - 6 \kappa_1 \kappa_3 + 4 (\kappa_2^2 + \kappa_2 \kappa_3 + 4 \kappa_3^2) \geq 0, \]
we have
\begin{equation} \label{eq6-3}
W_3 \geq W_3 \Big|_{\kappa_1 = \kappa_2} = (\kappa_2 - \kappa_3) (19 \kappa_2^2 + 4 \kappa_2 \kappa_3 + 4 \kappa_3^2) \geq 0.
\end{equation}
Now we require the following condition:
\begin{equation} \label{Condition7}
\frac{1}{\beta - 1} - \frac{1}{\beta} < \frac{1 - \varphi}{1 + \varphi}.
\end{equation}
By \eqref{eq6-4}, \eqref{eq6-3} and  \eqref{Condition7}, we arrive at
\[ T_3 \geq 0. \]

{\bf Subcase (2).}
When $\kappa_3 \leq 0$, we immediately obtain
\[ W_3 =  15 \kappa_1^3 - 3 \kappa_1^2 \kappa_3  - 4 \kappa_3 (2 \kappa_2 + \kappa_3)^2 + 4 \kappa_1 (\kappa_2^2 + \kappa_2 \kappa_3 + 4 \kappa_3^2) \geq 15 \kappa_1^3. \]
When $\kappa_3 \geq 0$, we have
\[ 0 \leq \kappa_3 \leq \kappa_2 \leq \phi \kappa_1.  \]
Then
\[ \begin{aligned}
W_3 \geq & 15 \kappa_1^3 - 3 \kappa_1^2 \phi \kappa_1  - 4 \phi \kappa_1 (2 \phi \kappa_1 + \phi \kappa_1)^2 \\
= & 3 \kappa_1^3 (5 - \phi - 12 \phi^3).
\end{aligned} \]
Now we impose the following condition:
\begin{equation} \label{Condition8}
\phi < \phi_0,
\end{equation}
where $\phi_0$ is the single real root of the cubic polynomial
\[ 5 - \phi - 12 \phi^3. \]
Mathematica shows that
\[ \phi_0 \approx 0.709742. \]
Then by \eqref{eq6-2} and \eqref{Condition7},
\begin{equation*}
\begin{aligned}
T_3 \geq & (\kappa_1 - \phi \kappa_1) \kappa_1 \Big( (\beta^2 - \beta - 1) \kappa_1 - (\beta^2 - \beta + 1) \varphi \kappa_1 \Big)  3 \kappa_1^3 (5 - \phi - 12 \phi^3) \\
& - C  \kappa_1^{- \frac{1}{3}} \kappa_1^5 \geq 0
\end{aligned}
\end{equation*}
when $\kappa_1$ is sufficiently large. Here we have applied Proposition \ref{Prop0}.

\vspace{2mm}

\subsection{The case when $\kappa_3 \leq - \varphi \kappa_1$.}~

\vspace{2mm}

Since
\[ \frac{2 \beta}{\nu^{4 + 1}} H  f_3 \frac{u_3^2}{u^2} \big( \kappa_3 - \nu^{4 + 1} \big) \]
can be eaten by $( \beta - 1 ) H f_3 \kappa_3^2$,
we want to guarantee the nonnegativity of the following quantity:
\begin{equation*}
\begin{aligned}
& E_3 := Q_3 - \frac{\beta (\beta - 1)}{(\nu^{4 + 1})^2} H  f_3  \frac{u_3^2}{u^2} (\nu^{4 + 1} - \kappa_3)^2 \\
= &  \frac{\sigma}{4} \frac{\beta}{(\nu^{4 + 1})^2} \frac{u_3^2}{u^2} (\nu^{4 + 1} - \kappa_3)^2 H \frac{S_3}{(\kappa_1 + \kappa_2 + \kappa_3) (\kappa_1 + \kappa_3 + \kappa_4) (\kappa_2 + \kappa_3 + \kappa_4) B_3}.
\end{aligned}
\end{equation*}
Now we want to prove the nonnegativity of $S_3$, which is a homogeneous polynomial of degree $6$.
Similar as the above argument, we can differentiate $S_3$ several times with respect to $\kappa_4$ to prove that
\[ \frac{\partial^2 S_3}{\partial \kappa_4^2} \geq 0. \]
Hence
\[ \begin{aligned}
& \frac{\partial S_3}{\partial \kappa_4} \geq \frac{\partial S_3}{\partial \kappa_4} \Big|_{\kappa_4 = - \kappa_2 - \kappa_3} := X_3.
\end{aligned} \]
By differentiating $X_3$ several times with respect to $\kappa_1$ we can prove that

\[ \begin{aligned}
& X_3 \geq X_3 \Big|_{\kappa_1 = \kappa_2} \\
= & 2 (\kappa_2 - \kappa_3) (2 \kappa_2 + \kappa_3) \Big( (1 + \beta) \kappa_2 + (-1 + \beta) \kappa_3 \Big) \Big( 19 \kappa_2^2 + 4 \kappa_2 \kappa_3 + 4 \kappa_3^2 \Big) \geq 0.
\end{aligned} \]
Also in view of \eqref{eq6-3}, we can deduce that
\[ S_3 \geq S_3 \Big|_{\kappa_4 = - \kappa_2 - \kappa_3} = (\kappa_1 - \kappa_2) (\kappa_1 + \kappa_2 + \kappa_3) \Big( (1 + \beta) \kappa_1 + (-1 + \beta) \kappa_3 \Big)  W_3 \geq 0.  \]

\vspace{4mm}

\section{For $i = 4$}

\vspace{4mm}

In this section, we divide our discussion into two cases. That is, the case when $\kappa_4 \geq - \theta \kappa_1$ and the case when $\kappa_4 \leq - \theta \kappa_1$. Here $0 < \theta < 2$ is a constant to be determined later.

\vspace{2mm}

\subsection{The case when $\kappa_4 \geq - \theta \kappa_1$.}~

\vspace{2mm}

By Cauchy-Schwartz inequality
\[  \frac{2 \beta}{\nu^{4 + 1}} H f_4 \frac{u_4^2}{u^2} \big( \kappa_4 - \nu^{4 + 1} \big) \geq - \frac{\beta^2}{(\beta - 1)(\nu^{4 + 1})^2} H f_4 \frac{u_4^2}{u^2} (\kappa_4 - \nu^{4 + 1})^2 - (\beta - 1) f_4 \frac{u_4^2}{u^2} H,   \]
we want to guarantee the nonnegativity of the following quantity:

\begin{equation*}
\begin{aligned}
& Q_4 - \frac{\beta^2}{(\beta - 1)(\nu^{4 + 1})^2} H f_4 \frac{u_4^2}{u^2} (\kappa_4 - \nu^{4 + 1})^2 - \frac{\beta (\beta - 1)}{(\nu^{4 + 1})^2} H f_4 \frac{u_4^2}{u^2} (\kappa_4 - \nu^{4 + 1})^2 \\
= &  \frac{\sigma}{4} \frac{\beta}{(\beta - 1) (\nu^{4 + 1})^2} \frac{u_4^2}{u^2} (\kappa_4 - \nu^{4 + 1})^2 H \frac{T_4}{ (\kappa_1 + \kappa_2 + \kappa_4) (\kappa_1 + \kappa_3 + \kappa_4) (\kappa_2 + \kappa_3 + \kappa_4) B_4},
\end{aligned}
\end{equation*}
where
\begin{equation*}
\begin{aligned}
B_4 := & 15 \kappa_1^4 + 15 \kappa_2^4 + 15 \kappa_3^4 + 2 \kappa_2^3 (8 \kappa_3 - \kappa_4) +
   2 \kappa_1^3 (8 \kappa_2 + 8 \kappa_3 - \kappa_4) \\
   & - 2 \kappa_3^3 \kappa_4 - 5 \kappa_3^2 \kappa_4^2 + 20 \kappa_3 \kappa_4^3 +
   12 \kappa_4^4 + \kappa_2^2 (2 \kappa_3^2 - 22 \kappa_3 \kappa_4 - 5 \kappa_4^2) \\
   & + \kappa_1^2 (2 \kappa_2^2 + 2 \kappa_3^2 - 22 \kappa_2 \kappa_4 - 22 \kappa_3 \kappa_4 - 5 \kappa_4^2) + 2 \kappa_2 (8 \kappa_3^3 - 11 \kappa_3^2 \kappa_4 + \kappa_3 \kappa_4^2 + 10 \kappa_4^3) \\
   & + 2 \kappa_1 \Big( 8 \kappa_2^3 + 8 \kappa_3^3 - 11 \kappa_2^2 \kappa_4 - 11 \kappa_3^2 \kappa_4 + \kappa_3 \kappa_4^2 + 10 \kappa_4^3 + \kappa_2 \kappa_4 (-36 \kappa_3 + \kappa_4) \Big),
\end{aligned}
\end{equation*}
\begin{equation*}
\begin{aligned}
T_4 := & (\kappa_1 + \kappa_3 + \kappa_4) (\kappa_1 + \kappa_2 +
    \kappa_4) \Big( (\beta^2 - \beta - 1) \kappa_1 + (\beta^2 - \beta + 1) \kappa_4 \Big) \\
   & \cdot \Big( 15 \kappa_1^3 +
    3 \kappa_1^2 (\kappa_2 + \kappa_3) + 4 (\kappa_2 - \kappa_3)^2 (\kappa_2 + \kappa_3) +
    4 \kappa_1 (4 \kappa_2^2 + 7 \kappa_2 \kappa_3 + 4 \kappa_3^2) \Big) \\
    & + (\kappa_2 + \kappa_3 + \kappa_4) M_4,
\end{aligned}
\end{equation*}
and $M_4$ is a homogeneous polynomial of degree $5$.

In order to make $T_4$ to be nonnegative, we impose the following conditions:
\begin{equation} \label{Condition3}
 0 < \theta < \frac{1}{2}
\end{equation}
and
\begin{equation} \label{Condition4}
\frac{1}{\beta - 1} - \frac{1}{\beta} < \frac{1 - \theta}{1 + \theta}.
\end{equation}
Consequently, we have
\begin{equation*}
\begin{aligned}
T_4 \geq & (\kappa_1 - \theta \kappa_1 - \theta \kappa_1) (\kappa_1 - \theta
    \kappa_1) \Big( (\beta^2 - \beta - 1) \kappa_1 - (\beta^2 - \beta + 1) \theta \kappa_1 \Big) 15 \kappa_1^3 \\
    & - C  \kappa_1^{- \frac{1}{3}} \kappa_1^5 \geq 0
\end{aligned}
\end{equation*}
when $\kappa_1$ is sufficiently large. Here we have applied Proposition \ref{Prop0}.

\vspace{2mm}

\subsection{The case when $\kappa_4 \leq - \theta \kappa_1$.}

\vspace{2mm}

\begin{equation*}
\begin{aligned}
Q_4 = & - \sum\limits_{j l} f_{jl}  h_{jj4} h_{ll4} +  \frac{2 \sigma}{4} \frac{1}{H - \kappa_4}  \frac{1}{H - \kappa_1} t_1^2 + \frac{2 \sigma}{4} \frac{1}{H - \kappa_4} \frac{1}{H - \kappa_2} t_2^2 \\
& + \frac{2 \sigma}{4} \frac{1}{H - \kappa_4} \frac{1}{H - \kappa_3} t_3^2 \\
= &  \frac{\sigma}{4} \frac{\beta^2}{(\nu^{4 + 1})^2} \frac{u_4^2}{u^2} (\nu^{4 + 1} - \kappa_4)^2 H^2 \\
& \cdot \frac{\left( \substack{ 2 (\kappa_2 + \kappa_4 + \kappa_3) (\kappa_1 + \kappa_4 + \kappa_3) A_3 +  2 (\kappa_2 + \kappa_4 + \kappa_3) (\kappa_1 + \kappa_4 + \kappa_2) A_2 \\
+ 2 (\kappa_2 + \kappa_4 + \kappa_1) (\kappa_1 + \kappa_4 + \kappa_3) A_1 +  4 (\kappa_2 + \kappa_4 + \kappa_3) (\kappa_1 + \kappa_4 + \kappa_3) (\kappa_1 + \kappa_4 + \kappa_2) A_4 } \right)}{(\kappa_1 + \kappa_2 + \kappa_4) (\kappa_1 + \kappa_3 + \kappa_4) (\kappa_2 + \kappa_3 +
     \kappa_4) B_4 },
\end{aligned}
\end{equation*}
where
\begin{equation*}
\begin{aligned}
A_3 = & 4 \kappa_1^3 - 4 \kappa_1^2 \kappa_2 - 4 \kappa_1 \kappa_2^2 +
    4 \kappa_2^3 + 16 \kappa_1^2 \kappa_3 + 28 \kappa_1 \kappa_2 \kappa_3 \\
    & + 16 \kappa_2^2 \kappa_3 + 3 \kappa_1 \kappa_3^2 +
    3 \kappa_2 \kappa_3^2 + 15 \kappa_3^3,
\end{aligned}
\end{equation*}

\begin{equation*}
\begin{aligned}
A_2 = & 4 \kappa_1^3 + 16 \kappa_1^2 \kappa_2 + 3 \kappa_1 \kappa_2^2 +
    15 \kappa_2^3 - 4 \kappa_1^2 \kappa_3 + 28 \kappa_1 \kappa_2 \kappa_3 \\
    & + 3 \kappa_2^2 \kappa_3 - 4 \kappa_1 \kappa_3^2 + 16 \kappa_2 \kappa_3^2 + 4 \kappa_3^3,
\end{aligned}
\end{equation*}

\begin{equation*}
\begin{aligned}
A_1 = & 15 \kappa_1^3 + 3 \kappa_1^2 \kappa_2 + 16 \kappa_1 \kappa_2^2 +
    4 \kappa_2^3 + 3 \kappa_1^2 \kappa_3 + 28 \kappa_1 \kappa_2 \kappa_3 \\
    & - 4 \kappa_2^2 \kappa_3 + 16 \kappa_1 \kappa_3^2 -
    4 \kappa_2 \kappa_3^2 + 4 \kappa_3^3,
\end{aligned}
\end{equation*}

\begin{equation*}
\begin{aligned}
A_4 = & 85 \kappa_1^2 + 50 \kappa_1 \kappa_2 + 19 \kappa_2^2 + 50 \kappa_1 \kappa_3 + 20 \kappa_2 \kappa_3 + 19 \kappa_3^2 \\
& + (\kappa_2 + \kappa_3 + \kappa_4)(102 \kappa_1 + 66 \kappa_2 + 66 \kappa_3 + 36 \kappa_4),
\end{aligned}
\end{equation*}

\begin{equation*}
\begin{aligned}
B_4
= & (\kappa_1 + \kappa_2 + \kappa_3) A_1 + (\kappa_2 + \kappa_3 + \kappa_4) P_4,
\end{aligned}
\end{equation*}
where $P_4$ is a homogeneous polynomial of degree $3$.

Since
\[ \frac{2 \beta}{\nu^{4 + 1}} H  f_4 \frac{u_4^2}{u^2} \big( \kappa_4 - \nu^{4 + 1} \big) \]
can be eaten by $( \beta - 1 ) H f_4 \kappa_4^2$,
we want to guarantee the nonnegativity of the following quantity:
\begin{equation*}
\begin{aligned}
& E_4 := Q_4 - \frac{\beta (\beta - 1)}{(\nu^{4 + 1})^2} H  f_4  \frac{u_4^2}{u^2} (\nu^{4 + 1} - \kappa_4)^2 \\
= &  \frac{\sigma}{4} \frac{\beta^2}{(\nu^{4 + 1})^2} \frac{u_4^2}{u^2} (\nu^{4 + 1} - \kappa_4)^2 H^2  \Bigg(  \frac{1}{(\kappa_1 + \kappa_2 + \kappa_4)} \bigg( \frac{2  A_3}{ B_4} - \Big( 1 - \frac{1}{\beta} \Big) \frac{1}{H} \bigg) \\
& +  \frac{1}{(\kappa_1 + \kappa_3 + \kappa_4)} \bigg( \frac{2 A_2}{B_4} - \Big( 1 - \frac{1}{\beta} \Big) \frac{1}{H} \bigg) \\
& + \frac{1}{(\kappa_2 + \kappa_3 +
     \kappa_4)} \bigg( \frac{2 A_1}{ B_4} - \Big( 1 - \frac{1}{\beta} \Big) \frac{1}{H} \bigg) + \frac{4 A_4}{B_4} \Bigg),
\end{aligned}
\end{equation*}
which can be rewritten as
\begin{equation*}
\begin{aligned}
& E_4 = \frac{\sigma}{4} \frac{\beta^2}{(\nu^{4 + 1})^2} \frac{u_4^2}{u^2} (\nu^{4 + 1} - \kappa_4)^2 H^2  \\ & \cdot \Bigg(  \frac{1}{(\kappa_1 + \kappa_2 + \kappa_4)} \bigg( \frac{2  A_3 + 2 A_4 (\kappa_1 - \kappa_3) + (\kappa_2 + \kappa_3 + \kappa_4) 2 A_4}{ B_4} - \Big( 1 - \frac{1}{\beta} \Big) \frac{1}{H} \bigg) \\
& +  \frac{1}{(\kappa_1 + \kappa_3 + \kappa_4)} \bigg( \frac{2 A_2 + 2 A_4 (\kappa_1 - \kappa_2) + (\kappa_2 + \kappa_3 + \kappa_4) 2 A_4}{B_4} - \Big( 1 - \frac{1}{\beta} \Big) \frac{1}{H} \bigg) \\
& + \frac{1}{(\kappa_2 + \kappa_3 +
     \kappa_4)} \bigg( \frac{2 A_1}{ B_4} - \Big( 1 - \frac{1}{\beta} \Big) \frac{1}{H} \bigg) \Bigg) .
\end{aligned}
\end{equation*}
We note that
\begin{equation*}
\begin{aligned}
A_4 = A_5 + (\kappa_2 + \kappa_3 + \kappa_4)(102 \kappa_1 + 66 \kappa_2 + 66 \kappa_3 + 36 \kappa_4),
\end{aligned}
\end{equation*}
where
\begin{equation*}
\begin{aligned}
A_5 :=  85 \kappa_1^2 + 50 \kappa_1 \kappa_2 + 19 \kappa_2^2 + 50 \kappa_1 \kappa_3 + 20 \kappa_2 \kappa_3 + 19 \kappa_3^2.
\end{aligned}
\end{equation*}
Hence $E_4$ can be rewritten as
\begin{equation} \label{eq5-6}
\begin{aligned}
& E_4 = \frac{\sigma}{4} \frac{\beta^2}{(\nu^{4 + 1})^2} \frac{u_4^2}{u^2} (\nu^{4 + 1} - \kappa_4)^2 H^2  \\ & \cdot \Bigg(  \frac{1}{(\kappa_1 + \kappa_2 + \kappa_4)} \bigg( \frac{2  A_3 + 2 A_5 (\kappa_1 - \kappa_3) + (\kappa_2 + \kappa_3 + \kappa_4) P_5}{ B_4} - \Big( 1 - \frac{1}{\beta} \Big) \frac{1}{H} \bigg) \\
& +  \frac{1}{(\kappa_1 + \kappa_3 + \kappa_4)} \bigg( \frac{2 A_2 + 2 A_5 (\kappa_1 - \kappa_2) + (\kappa_2 + \kappa_3 + \kappa_4) P_5}{B_4} - \Big( 1 - \frac{1}{\beta} \Big) \frac{1}{H} \bigg) \\
& + \frac{1}{(\kappa_2 + \kappa_3 +
     \kappa_4)} \bigg( \frac{2 A_1}{ B_4} - \Big( 1 - \frac{1}{\beta} \Big) \frac{1}{H} \bigg) \Bigg),
\end{aligned}
\end{equation}
where $P_5$ is a homogeneous polynomial of degree $2$ which may change from line to line.

\begin{lemma} \label{estimate of B4}
We have\begin{enumerate}
         \item  \[ B_4 \leq  (243 + \delta) \kappa_1^4 \]
when $\kappa_1$ is sufficiently large. Here $\delta > 0$ is the same with the one in Condition \eqref{Condition1};

         \item  \[  81 \kappa_1^3 \geq A_1 \geq 3 (5 + \theta + 5 \theta^2) \kappa_1^3. \]
       \end{enumerate}
\end{lemma}
\begin{proof}
Consider $A_1$ as a polynomial of $\kappa_2$ and $\kappa_3$:
\begin{equation*}
\begin{aligned}
A_1 = & 15 \kappa_1^3 + 3 \kappa_1^2 \kappa_2 + 16 \kappa_1 \kappa_2^2 + 4 \kappa_2^3 + 3 \kappa_1^2 \kappa_3 \\
& + 28 \kappa_1 \kappa_2 \kappa_3 - 4 \kappa_2^2 \kappa_3 + 16 \kappa_1 \kappa_3^2 - 4 \kappa_2 \kappa_3^2 + 4 \kappa_3^3.
\end{aligned}
\end{equation*}
We want to find the extreme values of $A_1$ in the bounded closed set
\[ \mathcal{D}_4 := \Big\{ (\kappa_2, \kappa_3) \, \Big| \, \frac{\theta}{2} \kappa_1 \leq \kappa_2 \leq \kappa_1, \ \  \theta \kappa_1 - \kappa_2 \leq \kappa_3 \leq \kappa_2 \Big\}. \]
By elementary calculation,  $A_1$ does not have critical points in the interior of $\mathcal{D}_4$. Hence we find the extreme value of $A_1$ on the boundary of $\mathcal{D}_4$.

We first set
\[ \kappa_2 = \kappa_1, \quad (\theta - 1) \kappa_1 \leq \kappa_3 \leq \kappa_1 \]
and find
$A_1 \Big|_{\kappa_2 = \kappa_1}$, which is a function of $\kappa_3$, does not have any critical point on the open interval
\[ \Big( (\theta - 1) \kappa_1, \kappa_1 \Big). \]

We then set
\[ \kappa_3 = \kappa_2, \quad  \frac{\theta}{2} \kappa_1 \leq \kappa_2 \leq \kappa_1, \]
and find
\begin{equation*}
\begin{aligned}
A_1 \Big|_{\kappa_3 = \kappa_2} & = 3 \kappa_1 (5 \kappa_1^2 + 2 \kappa_1 \kappa_2 + 20 \kappa_2^2),
\end{aligned}
\end{equation*}
as a function of $\kappa_2$, does not have any critical point on the open interval $\Big(  \frac{\theta}{2} \kappa_1, \kappa_1 \Big)$.

Next, we set
\[ \kappa_3 =  \theta \kappa_1 - \kappa_2, \quad  \frac{\theta}{2} \kappa_1 \leq \kappa_2 \leq \kappa_1, \]
and find that
$A_1 \Big|_{\kappa_3 =  \theta \kappa_1 - \kappa_2}$,
as a function of $\kappa_2$, does not have a critical point on the open interval $\Big(  \frac{\theta}{2} \kappa_1, \kappa_1 \Big)$.

Finally, we compare the three values:

\[ A_1 \Big|_{\kappa_2 = \kappa_3 = \frac{\theta}{2} \kappa_1} = 3 (5 + \theta + 5 \theta^2) \kappa_1^3, \]

\[ A_1 \Big|_{\kappa_2 = \kappa_1, \, \kappa_3 = (\theta - 1) \kappa_1} = (1 + \theta) (19 - 4 \theta + 4 \theta^2) \kappa_1^3,  \]

\[ A_1 \Big|_{\kappa_2 = \kappa_1, \, \kappa_3 = \kappa_1} = 81 \kappa_1^3. \]

\end{proof}

By Lemma \ref{estimate of B4} (2), we can deduce that
\begin{equation}  \label{First part}
\frac{2 A_1}{ B_4} - \Big( 1 - \frac{1}{\beta} \Big) \frac{1}{H} \geq \bigg( \frac{2}{3} - \delta - 1 + \frac{1}{\beta} \bigg) \frac{1}{\kappa_1}.
\end{equation}
Hence we require that
\begin{equation} \label{Condition2}
\frac{2}{3} - \delta - 1 + \frac{1}{\beta} > 0.
\end{equation}
For this, we require that
\begin{equation} \label{Condition2-1}
\beta < 3.
\end{equation}

Similarly, we can prove that
\begin{lemma} \label{Estimate of second}
\[  \frac{2 A_2 + 2 A_5 (\kappa_1 - \kappa_2)}{(\kappa_1 + \kappa_2 + \kappa_3) A_1} \geq  \frac{2}{3 \kappa_1}. \]
\end{lemma}
\begin{proof}
We find the minimum value of
\[  \frac{2 A_2 + 2 A_5 (\kappa_1 - \kappa_2)}{(\kappa_1 + \kappa_2 + \kappa_3) A_1} \]
in the bounded closed set
\[ \mathcal{D}_4 = \Big\{ (\kappa_2, \kappa_3) \, \Big| \, \frac{\theta}{2} \kappa_1 \leq \kappa_2 \leq \kappa_1, \ \  \theta \kappa_1 - \kappa_2 \leq \kappa_3 \leq \kappa_2 \Big\}. \]
By elementary calculation, we know that there is no critical point in the interior of $\mathcal{D}_4$, neither does there exist critical point in the interior of the three sides of $\mathcal{D}_4$. We thus compare the values on the three vertices of $\mathcal{D}_4$:

\[  \frac{2 A_2 + 2 A_5 (\kappa_1 - \kappa_2)}{(\kappa_1 + \kappa_2 + \kappa_3) A_1} \Big|_{\kappa_2 = \kappa_3 = \frac{\theta}{2} \kappa_1} = \frac{356 + 54 \theta - 15 \theta^2 - 10 \theta^3}{6 (1 + \theta) (5 + \theta + 5 \theta^2) \kappa_1}, \]

\[  \frac{2 A_2 + 2 A_5 (\kappa_1 - \kappa_2)}{(\kappa_1 + \kappa_2 + \kappa_3) A_1} \Big|_{\kappa_2 = \kappa_1, \, \kappa_3 = (\theta - 1) \kappa_1} = \frac{2}{(1 + \theta) \kappa_1},  \]

\[  \frac{2 A_2 + 2 A_5 (\kappa_1 - \kappa_2)}{(\kappa_1 + \kappa_2 + \kappa_3) A_1} \Big|_{\kappa_2 = \kappa_1, \, \kappa_3 = \kappa_1} = \frac{2}{3 \kappa_1}. \]

\end{proof}

By Lemma \ref{Estimate of second} and Lemma \ref{estimate of B4} (2), we can deduce that
\begin{equation} \label{Second part}
\begin{aligned}
&  \frac{2 A_2 + 2 A_5 (\kappa_1 - \kappa_2) + (\kappa_2 + \kappa_3 + \kappa_4) P_5}{B_4} - \Big( 1 - \frac{1}{\beta} \Big) \frac{1}{H} \\
 \geq & \bigg( \frac{2}{3} - \delta - 1 + \frac{1}{\beta} \bigg) \frac{1}{\kappa_1}.
\end{aligned}
\end{equation}

\begin{lemma} \label{Estimate of third}
\[  \frac{2  A_3 + 2 A_5 (\kappa_1 - \kappa_3)}{ (\kappa_1 + \kappa_2 + \kappa_3) A_1} \geq  \frac{2}{3 \kappa_1}. \]
\end{lemma}
\begin{proof}

We find the minimum value of
\[ \frac{2  A_3 + 2 A_5 (\kappa_1 - \kappa_3)}{(\kappa_1 + \kappa_2 + \kappa_3) A_1} \]
in the bounded closed set
\[ \mathcal{D}_4 = \Big\{ (\kappa_2, \kappa_3) \, \Big| \, \frac{\theta}{2} \kappa_1 \leq \kappa_2 \leq \kappa_1, \ \  \theta \kappa_1 - \kappa_2 \leq \kappa_3 \leq \kappa_2 \Big\}. \]
By elementary calculation, we know that there is no critical point in the interior of $\mathcal{D}_4$. On the side
\[ \Big\{ (\kappa_2, \kappa_3) \Big| \kappa_2 + \kappa_3 = \theta \kappa_1, \quad  \frac{\theta}{2} \kappa_1 \leq \kappa_2 \leq \kappa_1 \Big\}, \]
for the function
\[ \frac{2  A_3 + 2 A_5 (\kappa_1 - \kappa_3)}{(\kappa_1 + \kappa_2 + \kappa_3) A_1} \Big|_{\kappa_3 = \theta \kappa_1 - \kappa_2}, \]
there is a single local maximum point in $(\frac{\theta}{2} \kappa_1, \kappa_1)$. Except this, there is no other critical point in the interior of the three sides of $\mathcal{D}_4$. We thus compare the values on the three vertices of $\mathcal{D}_4$:

\[ \frac{2  A_3 + 2 A_5 (\kappa_1 - \kappa_3)}{(\kappa_1 + \kappa_2 + \kappa_3) A_1} \Big|_{\kappa_2 = \kappa_3 = \frac{\theta}{2} \kappa_1} = \frac{356 + 54 \theta - 15 \theta^2 - 10 \theta^3}{6 (1 + \theta) (5 + \theta + 5 \theta^2) \kappa_1}, \]

\[ \frac{2  A_3 + 2 A_5 (\kappa_1 - \kappa_3)}{(\kappa_1 + \kappa_2 + \kappa_3) A_1} \Big|_{\kappa_2 = \kappa_1, \, \kappa_3 = (\theta - 1) \kappa_1} = \frac{274 + 108 \theta - 66 \theta^2 - 8 \theta^3}{(1 + \theta)^2 (19 - 4 \theta + 4 \theta^2) \kappa_1},  \]

\[ \frac{2  A_3 + 2 A_5 (\kappa_1 - \kappa_3)}{(\kappa_1 + \kappa_2 + \kappa_3) A_1} \Big|_{\kappa_2 = \kappa_1, \, \kappa_3 = \kappa_1} = \frac{2}{3 \kappa_1}. \]

\end{proof}

By Lemma \ref{Estimate of third} and Lemma \ref{estimate of B4} (2), we can deduce that
\begin{equation} \label{Third part}
\begin{aligned}
& \frac{2  A_3 + 2 A_5 (\kappa_1 - \kappa_3) + (\kappa_2 + \kappa_3 + \kappa_4) P_5}{ B_4} - \Big( 1 - \frac{1}{\beta} \Big) \frac{1}{H} \\
\geq &  \bigg( \frac{2}{3} - \delta - 1 + \frac{1}{\beta} \bigg) \frac{1}{\kappa_1}.
\end{aligned}
\end{equation}

In view of \eqref{eq5-6}, \eqref{First part}, \eqref{Second part} and \eqref{Third part} and requiring \eqref{Condition2}, we have
\begin{equation*}
\begin{aligned}
& E_4 \geq 0.
\end{aligned}
\end{equation*}

\vspace{4mm}

\section{The final step}

\vspace{4mm}

Now we choose
\[ \beta = 2.999, \quad \theta = 0.499, \quad \phi = 0.7, \quad \varphi = 0.4.  \]
We can verify that Condition \eqref{Condition1-1}, \eqref{Condition2-1}, \eqref{Condition3}, \eqref{Condition4}, \eqref{Condition5}, \eqref{Condition6}, \eqref{Condition7}, \eqref{Condition8} are satisfied. Then we choose $\delta > 0$ sufficiently small such that Condition \eqref{Condition1}, \eqref{Condition2} are satisfied. Then we obtain a uniform upper bound for $\kappa_1$ by \eqref{eq5-2}.

\end{document}